\documentclass[12pt]{article}

\usepackage{amsmath,amssymb,amsthm,amssymb}
\usepackage{tikz,tkz-graph,tikz-cd,tikz-qtree}
\usepackage{youngtab}
\usepackage[OT2,T1]{fontenc}
\usepackage{indentfirst}

\DeclareSymbolFont{cyrletters}{OT2}{wncyr}{m}{n}

\DeclareMathSymbol{\Sha}{\mathalpha}{cyrletters}{"58}

\newtheorem{theorem}{Theorem}[section]
\newtheorem{lemma}{Lemma}[section]
\newtheorem{corollary}{Corollary}[section]

\newtheorem{definition}{Definition}[section]
\newtheorem{proposition}{Proposition}[section]
\newtheorem{remark}{Remark}[section]

\newtheorem{examples}{Example}[section]

\begin{document}
	
\bibliographystyle{abbrv}
	
\title{Explicit Representatives and Sizes of Cyclotomic Cosets \uppercase\expandafter{\romannumeral2}: 
Cyclotomic Systems and Applications to Constacyclic Codes}
	
\author{Li Zhu$^{1}$ and Hongfeng Wu$^{2}$\footnote{Corresponding author.}
\setcounter{footnote}{-1}
\footnote{E-mail addresses: lizhumath@pku.edu.cn (L. Zhu),
whfmath@gmail.com (H. Wu).}
\\[0.5ex]
\small $^{1}$School of Mathematical Sciences, Guizhou Normal University, Guiyang, China
\\
\small $^{2}$College of Science, North China University of Technology, Beijing, China}
	
\date{}
\maketitle
	
\thispagestyle{plain}
\setcounter{page}{1}

\begin{abstract}
Let $q=p^e$ be a prime power, let $\ell$ be a prime with 
$\gcd(\ell,qn)=1$, and let $n$ be a positive integer. While 
$q$-cyclotomic cosets are usually considered modulo a fixed integer, 
we consider their behavior along the sequence $n,\ell n,\ell^2n,\ldots$ of moduli. 
To describe the resulting compatibility among cyclotomic cosets, we 
introduce the $\ell$-adic $q$-cyclotomic system with base module $n$, 
defined as the projective limit of the spaces of $q$-cyclotomic cosets 
modulo $\ell^i n$. We associate to each $q$-cyclotomic coset modulo $n$ 
a cyclotomic $\ell$-adic integer and use its $\ell$-adic expansion to  
describe all compatible liftings through the successive levels. This 
leads to a classification of the corresponding sequences according to 
their splitting and stable behaviors, together with explicit formulas 
for the representatives and sizes of their components. We treat 
separately the cases of odd $\ell$ and $\ell=2$. These results are then 
used to determine representatives and sizes of $q$-cyclotomic cosets 
for arbitrary admissible parameters. As an application, we obtain 
explicit irreducible factorizations of binomials over finite fields and 
use them to describe the corresponding constacyclic codes.

\end{abstract}

\section{Introduction}
Let $q$ be a prime power and let $n$ be a positive integer
coprime to $q$. The $q$-cyclotomic cosets modulo $n$ arise from the
action of the Frobenius map on the $n$-th roots of unity over
$\mathbb{F}_q$. In particular, they determine the irreducible
factorization of $X^n-1$ over $\mathbb{F}_q$, and hence enter naturally into the study of
cyclic and constacyclic codes. For this reason, explicit descriptions
of cyclotomic cosets, including their representatives and sizes, have
been considered in various settings.

Most existing descriptions of cyclotomic cosets concern a fixed
modulus. For particular families of parameters, many results on explicit
representatives and sizes can be obtained by using the arithmetic
properties of the modulus and the corresponding multiplicative
orders. In a series of papers, for instance \cite{Chen}, \cite{Sharma}, \cite{Liu_2}, \cite{Wu_2}, 
etc., the representatives and sizes of the corresponding classes of cyclotomic cosets are 
determined in order to calculate certain constacyclic codes. In \cite{Chen_3} 
and \cite{Wang} the $q$-cyclotomic cosets contained in the subset $1+r\mathbb{Z}/nr\mathbb{Z}$ 
of $\mathbb{Z}/nr\mathbb{Z}$, where $\mathrm{gcd}(q,nr)=1$ and $r \mid q-1$, 
are characterized and enumerated. Applying the results on cyclotomic cosets, \cite{Chen_3} 
gives the enumeration of Euclidean self-dual codes, and \cite{Wang} construct 
several classes of $p^{h}$-LCD MDS codes. And in \cite{Yue} and \cite{Geraci}, 
concerning with stream cipher $m$-sequences and problems in statistic 
physics respectively, algorithms to calculate the leader of cyclotomic cosets are given.

In our previous work \cite{Zhu_2}, we determined a complete set of
representatives and the sizes of all $q$-cyclotomic cosets modulo an
arbitrary admissible integer. This leads naturally to a different
question when the modulus is allowed to vary. Given a prime
$\ell$ with $\gcd(\ell,qn)=1$, consider the sequence $n,\ell n, \ell^2n,\ldots$ 
of moduli. The canonical projections between these residue rings induce maps
between the corresponding sets of $q$-cyclotomic cosets. One may then
ask how a cyclotomic coset at one level is related to its images and
preimages at the neighboring levels, and how its size changes along
the sequence.

To keep track of the compatibility among the cyclotomic cosets at
different levels, we introduce the $\ell$-adic $q$-cyclotomic system
with base module $n$,
\[
\mathcal{PC}_n^\ell
=
\varprojlim_{i \in \mathbb{N}} \mathcal{C}_{\ell^i n},
\]
where $\mathcal{C}_{\ell^i n}$ denotes the space of $q$-cyclotomic cosets modulo
$\ell^{i}n$. Thus an element of $\mathcal{PC}_n^\ell$ is a compatible 
sequence $(c_{\ell^i n})_{i\geq0}$ of cyclotomic cosets. The lifting 
problem for such sequences can be described in terms of
$\ell$-adic arithmetic. More precisely, for a $q$-cyclotomic coset
$c_n(\gamma)$ modulo $n$, we associate the cyclotomic $\ell$-adic
integer
\[
\phi(\gamma)=-\frac{\gamma}{n}\in\mathbb{Z}_\ell.
\]
The $\ell$-adic expansion of $\phi(\gamma)$ determines the possible
liftings of $c_n(\gamma)$ to the successive moduli $\ell^i n$. With this 
description, we classify the compatible sequences in
$\mathcal{PC}_n^\ell$ according to their splitting and stable behaviors. We
introduce the quasi-stable degree and the stable degree of a sequence
and determine the stage at which its components begin to remain
stable, together with the sizes of the components at the successive
levels. The cases of odd $\ell$ and $\ell=2$ are treated separately,
reflecting the different forms of the relevant $\ell$-adic valuation
formulas. We also consider the total $q$-cyclotomic system obtained by
letting the base moduli vary.

The preceding results give an explicit procedure for determining
representatives and sizes of the $q$-cyclotomic cosets contained in
\[
1+g\mathbb{Z}/ng\mathbb{Z},
\]
where $g$ is a positive divisor of $q-1$. The case $g=1$ recovers the
usual problem of determining all $q$-cyclotomic cosets modulo $n$.
For general $g$, the resulting coset partition is used to describe
the irreducible factorization of binomials
\[
X^N-\lambda
\]
over finite fields. We then apply the factorization to the
description and enumeration of $\lambda$-constacyclic codes, including
the repeated-root case.

The paper is organized as follows. Section 2 recalls the necessary
background on finite fields, cyclotomic cosets, and $\ell$-adic
integers. In Section 3, we introduce the $\ell$-adic and total
$q$-cyclotomic systems. Sections 4--6 develop the $\ell$-adic
description of compatible sequences and their splitting and stable
behaviors. Section 7 gives explicit representatives and sizes of
cyclotomic cosets in $1+g\mathbb{Z}/ng\mathbb{Z}$, with separate
treatments of odd primes and $\ell=2$. Section 8 applies these results
to the factorization of binomials and the classification of
constacyclic codes. Section 9 gives an example illustrating the
resulting procedure, and Appendix A collects the equal-difference
results used in the factorization arguments.

\section{Basic definitions and notations}
Throughout this paper, we fix a prime power $q=p^{e}$, where $p$ is a prime number and $e \in \mathbb{N}^{+}$. For notations that are dependent on $q$, the subscript $q$ is often omitted.

For integers $a \leq b$, we denote by $[a,b]$ the integer interval $\{a,a+1,\cdots,b\}$.

For $n \in \mathbb{N}^{+}$, $\mathbb{Z}/n\mathbb{Z}$ denotes the ring of residue classes modulo $n$. As there is a canonical surjection
$$\mathbb{Z} \rightarrow \mathbb{Z}/n\mathbb{Z}; \ \gamma \mapsto \overline{\gamma} = \gamma + n\mathbb{Z},$$
if there is no ambiguity, we will sometimes identify $\gamma$ with its image $\overline{\gamma}$ in $\mathbb{Z}/n\mathbb{Z}$. Therefore the notation $\gamma \in \mathbb{Z}/n\mathbb{Z}$ actually means that the image of $\gamma$ lies in $\mathbb{Z}/n\mathbb{Z}$.

If $n^{\prime}$ is a divisor of $n$, we denote by $\pi_{n/n^{\prime}}$ the natural projection
$$\pi_{n/n^{\prime}}: \mathbb{Z}/n\mathbb{Z} \rightarrow \mathbb{Z}/n^{\prime}\mathbb{Z}; \gamma+n\mathbb{Z} \mapsto \gamma + n^{\prime}\mathbb{Z}.$$

If $m$ and $n$ are coprime integers, we denote by $\mathrm{ord}_{n}(m)$ the multiplicative order of $m$ in $(\mathbb{Z}/ n\mathbb{Z})^{\ast}$, i.e., the smallest integer such that
$$m^{\mathrm{ord}_{n}(m)} \equiv 1 \pmod{n}.$$

Let $\ell$ be a prime number. We denote by $v_{\ell}(n)$ the $\ell$-adic valuation of $n$, i.e., the maximal integer such that $\ell^{v_{\ell}(n)} \mid n$. The following lift-the-exponent lemmas deal with the case where $\ell$ is odd and that $\ell = 2$ respectively.

\begin{lemma}{\cite{Nezami}}\label{lem 3}
	Let $\ell$ be an odd prime number, and $m$ be an integer such that $\ell \mid m-1$. Then $v_{\ell}(m^{d}-1) = v_{\ell}(m-1) + v_{\ell}(d)$ for any positive integer $d$.
\end{lemma}

\begin{lemma}{\cite{Nezami}}\label{lem 2}
	Let $m$ be an odd integer, and $d$ be a positive integer.
	\begin{itemize}
		\item[(1)] If $m \equiv 1 \pmod{4}$, then
		$$v_{2}(m^{d}-1) = v_{2}(m-1) + v_{2}(d), \  v_{2}(m^{d}+1) = 1.$$
		\item[(2)] If $m \equiv 3 \pmod{4}$ and $d$ is odd, then
		$$v_{2}(m^{d}-1) = 1, \  v_{2}(m^{d}+1) = v_{2}(m+1).$$
		\item[(3)] If $m \equiv 3 \pmod{4}$ and $d$ is even, then
		$$v_{2}(m^{d}-1) = v_{2}(m+1) + v_{2}(d), \  v_{2}(m^{d}+1) = 1.$$
	\end{itemize}
\end{lemma}

For a prime $\ell$, the ring $\mathbb{Z}_{\ell}$ of $\ell$-adic integers is the projective limit
$$\mathbb{Z}_{\ell} = \varprojlim_{i \in \mathbb{N}}\mathbb{Z}/\ell^{i}\mathbb{Z},$$
with respect to the canonical homomorphisms $\pi_{\ell^{i}/ \ell^{j}}: \mathbb{Z}/\ell^{i}\mathbb{Z} \rightarrow \mathbb{Z}/\ell^{j}\mathbb{Z}$, $i \geq j$. The projections $\mathbb{Z} \rightarrow \mathbb{Z}/\ell^{i}\mathbb{Z}$, $i \in \mathbb{N}$, induce an embedding $\mathbb{Z} \rightarrow \mathbb{Z}_{\ell}$, therefore $\mathbb{Z}$ can be regarded naturally as a subring of $\mathbb{Z}_{\ell}$.

The elements in $\mathbb{Z}_{\ell}$ are compatible sequences $a = (a_{i})_{i \in \mathbb{N}}$, where $a_{i} \in \mathbb{Z}/\ell^{i}\mathbb{Z}$ and where the compatibility means that for any $i \geq j$, $\pi_{\ell^{i}/ \ell^{j}}(a_{i}) = a_{j}$. Such sequences are in an one-to-one correspondence with $\ell$-adic power series via the map
$$a = (a_{i})_{i \in \mathbb{N}}\mapsto \sum_{i=0}^{\infty}a_{i}^{\prime}\cdot \ell^{i},$$
where for each $i \in \mathbb{N}$ the coefficient $a_{i}^{\prime}$ is the unique integer in the range $0 \leq a_{i}^{\prime} \leq \ell-1$ satisfying 
$$a_{i}^{\prime} \equiv \frac{a_{i+1}-a_{i}}{\ell^{i}} \pmod{\ell}.$$ 
Moreover, this correspondence preserves the additions and multiplications, that is, it gives a ring isomorphism from $\mathbb{Z}_{\ell}$ onto the ring of $\ell$-adic power series. The power series $\sum\limits_{i=0}^{\infty}a_{i}^{\prime}\cdot \ell^{i}$ is called the $\ell$-adic expansion of $a$. In this paper, we do not distinguish an $\ell$-adic integer from its $\ell$-adic expansion.

Let $\mathbb{F}_{q}$ be a finite field with $q$ elements, where $q$ is a prime power. The multiplicative group $\mathbb{F}_{q}^{\ast}$ of nonzero elements is a cyclic group. For any $\lambda \in \mathbb{F}_{q}^{\ast}$, the smallest positive integer $r$ such that $\lambda^{r} = 1$ is called the order of $\lambda$ and is denoted by $r = \mathrm{ord}(\lambda)$. It is obvious that $\mathrm{ord}(\lambda)$ is a divisor of $q-1$. 

Assume that $n$ is coprime to $q$. There are $n$ distinct roots of $X^{n}-1$ lying in some finite extension of $\mathbb{F}_{q}$, among which anyone that is not a root of $X^{m}-1$ for every $m < n$ is called a primitive $n$-th root of unity. In this note we fix a family
$$\{\zeta_{n} \ | \ \mathrm{gcd}(n,q)=1\},$$
where $\zeta_{n}$ is a primitive $n$-th root of unity, which is compatible in the sense that for any $m,n \in \mathbb{N}^{+}$ such that $m \mid n$ and $\mathrm{gcd}(n,q)=1$ it holds that 
$$\zeta_{n}^{\frac{n}{m}} = \zeta_{m}.$$
	
For any $\gamma \in \mathbb{Z}/ n\mathbb{Z}$, the $q$-cyclotomic coset modulo $n$ containing $\gamma$ is defined to be the subset
$$c_{n}(\gamma) = \{\gamma,\gamma q,\cdots,\gamma q^{\tau-1}\},$$
where $\tau$ is the smallest positive integer such that $\gamma q^{\tau} \equiv \gamma \pmod{n}$. Each element in $c_{n}(\gamma)$  is called a representative of the coset $c_{n}(\gamma)$. Regarding these elements as nonnegative integers less than $n$, the smallest one among them is called the leader of $c_{n}(\gamma)$. The order of $c_{n}(\gamma)$, which equals to $\tau$, is called the size of $c_{n}(\gamma)$. We denote by $\mathcal{C}_{n}$ the set of all $q$-cyclotomic cosets modulo $n$, and denote by $\mathcal{CR}_{n}$ a full set of representatives of $\mathcal{C}_{n}$.
	
It is a well-known fact that the $q$-cyclotomic cosets modulo $n$ are in an one-to-one correspondence with the irreducible factors of $X^{n}-1$ over $\mathbb{F}_{q}$. Concretely, to a $q$-cyclotomic coset $c_{n}(\gamma) = \{\gamma,\gamma q,\cdots,\gamma q^{\tau-1}\}$ we can assign an irreducible factor 
$$M_{\gamma}(X) = (X-\zeta_{n}^{\gamma})(X-\zeta_{n}^{\gamma q})\cdots(X-\zeta_{n}^{\gamma q^{\tau-1}})$$
of $X^{n}-1$ over $\mathbb{F}_{q}$. Then the irreducible factorization of $X^{n}-1$ over $\mathbb{F}_{q}$ is given by
$$X^{n}-1 = \prod_{\gamma \in \mathcal{CR}_{n}}M_{\gamma}(X).$$

Let $n$ be a positive integer, and $\lambda$ be a nonzero element in $\mathbb{F}_{q}$, with order $g = \mathrm{ord}(\lambda)$. A $\lambda$-constacyclic code of length $n$ over $\mathbb{F}_{q}$ is a $\mathbb{F}_{q}$-subspace of $\mathbb{F}_{q}^{n}$ that is closed under the $\lambda$-constacyclic shift map
$$s_{\lambda}: \mathbb{F}_{q}^{n} \rightarrow \mathbb{F}_{q}^{n}; \ (c_{0},c_{1},\cdots,c_{n-1}) \mapsto (\lambda c_{n-1},c_{0},\cdots,c_{n-2}).$$
Alternatively, a $\mathbb{F}$-subspace of $\mathbb{F}_{q}^{n}$ is $\lambda$-constacyclic if and only if it corresponds to an ideal in the quotient ring $\mathcal{R}_{n,\lambda} = \mathbb{F}_{q}[X]/(X^{n}-\lambda)$, under the bijection
$$\mathbb{F}_{q}^{n} \rightarrow \mathcal{R}_{n,\lambda}; \ (c_{0},c_{1},\cdots,c_{n-1}) \mapsto c_{0}+c_{1}X+\cdots+c_{n-1}X^{n-1}.$$
In this paper we will always identify a $\lambda$-constacyclic code $\mathcal{C}$ of length $n$ with an ideal $(f(X)) \subseteq \mathcal{R}_{n,\lambda}$, where $f(X)$ is a factor of $X^{n}-\lambda$. This polynomial $f(X)$ is called the generator polynomial of $\mathcal{C}$. Depending on whether or not $X^{n}-\lambda$ has repeated roots, a $\lambda$-constacyclic code $\mathcal{C}$ of length $n$ is said to be repeated-root if $p \mid n$, while is called simple-root otherwise.
	
\section{The definition of $q$-cyclotomic system}\label{sec 2}
Let $n$ and $n^{\prime}$ be positive integers such that $n^{\prime} \mid n$ and $\mathrm{gcd}(n,q) = 1$. The canonical projection $\pi_{n/n^{\prime}}: \mathbb{Z}/n\mathbb{Z} \rightarrow \mathbb{Z}/n^{\prime}\mathbb{Z}$ sends each $q$-cyclotomic coset modulo $n$ onto a $q$-cyclotomic coset modulo $n^{\prime}$, and consequently induces a map
$$\widehat{\pi}_{n/n^{\prime}}: \mathcal{C}_{n} \rightarrow \mathcal{C}_{n^{\prime}}; \ c_{n}(\gamma) \mapsto c_{n^{\prime}}(\pi_{n/n^{\prime}}(\gamma)) = c_{n^{\prime}}(\gamma).$$
Here the last equality holds as we identity $\pi_{n/n^{\prime}}(\gamma) = \gamma + n^{\prime}\mathbb{Z}$ with its preimage $\gamma$ in $\mathbb{Z}$.

Let $\ell$ be a prime number coprime to both $q$ and $n$. The sequence of projections
\begin{equation}\label{eq 6}
	\cdots \rightarrow \mathbb{Z}/\ell^{2}n\mathbb{Z} \rightarrow \mathbb{Z}/\ell n\mathbb{Z} \rightarrow \mathbb{Z}/n\mathbb{Z}
\end{equation}
induces a sequence of surjective maps
\begin{equation}\label{eq 1}
	\cdots \rightarrow \mathcal{C}_{\ell^{2}n} \rightarrow \mathcal{C}_{\ell n} \rightarrow \mathcal{C}_{n},
\end{equation}
which forms a projective system of finite sets and maps. Further, fixing a $q$-cyclotomic coset $c_{n}(\gamma)$ modulo $n$, the subsequence 
\begin{equation}\label{eq 7}
	\cdots \rightarrow \widehat{\pi}_{\ell^{2}n/n}^{-1}\{c_{n}(\gamma)\} \rightarrow \widehat{\pi}_{\ell n/n}^{-1}\{c_{n}(\gamma)\} \rightarrow \widehat{\pi}_{n/n}^{-1}\{c_{n}(\gamma)\} = \{c_{n}(\gamma)\}
\end{equation}
of \eqref{eq 1} also forms a projective system. Hence the following definitions are well-defined.

\begin{definition}
	Define the $\ell$-adic $q$-cyclotomic system with base module $n$ to be the projective limit
	$$\mathcal{PC}_{n}^{\ell} = \varprojlim_{i \in \mathbb{N}}\mathcal{C}_{\ell^{i}n}$$
	of the projective system \eqref{eq 1}. For a fixed $q$-cyclotomic coset $c_{n}(\gamma)$ modulo $n$, the $\ell$-adic $q$-cyclotomic system over $c_{n}(\gamma)$ is defined to be the projective limit
	$$\mathcal{PC}_{n}^{\ell}(\gamma) = \varprojlim_{i \in \mathbb{N}} \widehat{\pi}_{\ell^{i}n/n}^{-1}\{c_{n}(\gamma)\}$$
	of the projective system \eqref{eq 7}.
\end{definition}

Equivalently, the profinite space $\mathcal{PC}_{n}^{\ell}$ consists of compatible sequences $\mathfrak{c} = (c_{\ell^{i}n}(\gamma_{i}))_{i \in \mathbb{N}}$ of $q$-cyclotomic cosets, where the compatibility means that for any $i_{1} \geq i_{2}$ it holds that 
$$\pi_{\ell^{i_{1}}n/\ell^{i_{2}}n}(c_{\ell^{i_{1}}n}(\gamma_{i_{1}})) = c_{\ell^{i_{2}}n}(\gamma_{i_{2}}).$$
And $\mathcal{PC}_{n}^{\ell}(\gamma)$ is the sub-profinite space of $\mathcal{PC}_{n}^{\ell}$ consisting of the sequences with the first component being $c_{n}(\gamma)$. The following property is clear from the definition.

\begin{lemma}\label{lem 6}
	The $\ell$-adic $q$-cyclotomic system $\mathcal{PC}_{n}^{\ell}$ can be partitioned as 
	$$\mathcal{PC}_{n}^{\ell} = \bigsqcup_{\gamma \in \mathcal{CR}_{n}}\mathcal{PC}_{n}^{\ell}(\gamma),$$
	where $\mathcal{CR}_{n}$ is any full set of representatives of $q$-cyclotomic cosets modulo $n$.
\end{lemma}

Further, consider the projective system 
$$\{\mathcal{C}_{n} \ | \ \mathrm{gcd}(n,q)=1\},$$
equipped with the transition maps $\widehat{\pi}_{n_{1}/n_{2}}: \mathcal{C}_{n_{1}} \rightarrow \mathcal{C}_{n_{2}}$ where $n_{1}, n_{2}$ are such that $\mathrm{gcd}(n_{1},q)=1$ and $n_{2} \mid n_{1}$. We refer to the associated projective limit as the total $q$-cyclotomic system.

\begin{definition}
	The total $q$-cyclotomic system is defined to be the projective limit
	$$\mathcal{PC} = \varprojlim_{\mathrm{gcd}(n,q)=1} \mathcal{C}_{n}.$$
	Equivalently, $\mathcal{PC}$ is the profinite space consisting of all sequences $\mathfrak{c}=(c_{n}(\gamma_{n}))_{\mathrm{gcd}(n,q)=1}$ satisfying that 
	$$\pi_{n_{1}/n_{2}}(c_{n_{1}}(\gamma_{n_{1}})) = c_{n_{2}}(\gamma_{n_{2}})$$
	for all positive integers $n_{1}, n_{2}$ such that $\mathrm{gcd}(n_{1},q)=1$ and $n_{2} \mid n_{1}$.
\end{definition}

\section{Cyclotomic $\ell$-adic integer and cyclotomic $\ell$-adic class}
In this section we assume that 
\begin{itemize}
	\item[(1)] $q=p^{e}$ where $p$ is a prime number and $e$ is a positive integer;
	\item[(2)] $\ell$ is a prime number different than $p$; and
	\item[(3)] $n$ is a positive integer which is divisible by neither $p$ nor $\ell$.
\end{itemize}
For any $\overline{\gamma} = \gamma + n\mathbb{Z} \in \mathbb{Z}/n\mathbb{Z}$, we introduce the definitions of a cyclotomic $\ell$-adic integer over $\overline{\gamma}$ and of the cyclotomic $\ell$-adic class associated to $\overline{\gamma}$, which serves as a preparation for the classification of the sequences lying in $\mathcal{PC}_{n}^{\ell}(\gamma)$. 

Let us begin by determining the preimage $\pi_{\ell m/m}^{-1}(c_{m}(\gamma))$ of an arbitrary cyclotomic coset $c_{m}(\gamma)$ along $\pi_{\ell m/m}$, where $m$ is any positive integer coprime to $q$. First consider the case that $\ell$ divides $q^{\tau}-1$ for $\tau = |c_{m}(\gamma)|$.

\begin{lemma}\label{lem 4}
	Let $m$ be a positive integer coprime to $q$, and $c_{m}(\gamma)$ be a $q$-cyclotomic coset modulo $m$ of size $\tau$. Assume that $\ell$ divides $q^{\tau}-1$.
	\begin{itemize}
		\item[(1)] If $v_{\ell}(\gamma) + v_{\ell}(q^{\tau}-1) \geq v_{\ell}(m)+1$, then for any $0 \leq \sigma \leq \ell-1$ the $q$-cyclotomic coset modulo $\ell m$ containing $\gamma+\sigma m$ is given by
		$$c_{\ell m}(\gamma+\sigma m) = \{\gamma+\sigma m,(\gamma+\sigma m)q,\cdots,(\gamma+\sigma m)q^{\tau-1}\}.$$
		Moreover, the $c_{\ell m}(\gamma+\sigma m)$'s are pairwise distinct, and the preimage of $c_{m}(\gamma)$ along $\pi_{\ell m/m}$ is 
		\begin{equation}\label{eq 8}
			\pi_{\ell m/m}^{-1}(c_{m}(\gamma)) = \bigsqcup_{\sigma=0}^{\ell-1}c_{\ell m}(\gamma+\sigma m).
		\end{equation}
		\item[(2)] If $v_{\ell}(\gamma) + v_{\ell}(q^{\tau}-1) < v_{\ell}(m)+1$, then the $q$-cyclotomic coset modulo $\ell m$ containing $\gamma$ is given by
		$$c_{\ell m}(\gamma) = \{\gamma,\gamma q,\cdots,\gamma q^{\ell\tau-1}\}.$$
		Moreover, the preimage of $c_{m}(\gamma)$ along $\pi_{\ell m/m}$ is 
		$$\pi_{\ell m/m}^{-1}(c_{m}(\gamma)) = c_{\ell m}(\gamma) = c_{\ell m}(\gamma+m) = \cdots = c_{\ell m}(\gamma+(\ell-1)m).$$
	\end{itemize}
\end{lemma}

\begin{proof}
	\begin{itemize}
		\item[(1)] Since $\tau$ is the smallest positive integer such that $\gamma q^{\tau} \equiv \gamma \pmod{m}$ and $\ell m \mid \gamma q^{\tau}-\gamma$, then $\tau$ is also the smallest positive integer such that $\gamma q^{\tau} \equiv \gamma \pmod{\ell m}$. Therefore we have
		$$c_{\ell m}(\gamma) = \{\gamma,\gamma q,\cdots,\gamma q^{\tau-1}\}.$$
		For $1 \leq \sigma \leq \ell-1$, as $\ell \mid q^{\tau}-1$, it is clear that 
		$$(\gamma+\sigma m)q^{\tau} = \gamma q^{\tau}+\sigma mq^{\tau} \equiv \gamma + \sigma m \pmod{\ell m}.$$
		On the other hand, if $j$ is any positive integer satisfying the congruence
		\begin{equation}\label{eq 3}
			(\gamma+\sigma m)q^{j} \equiv \gamma+\sigma m \pmod{\ell m}
		\end{equation}
		then $\gamma(q^{j}-1) \equiv \sigma m(1-q^{j}) \pmod{\ell m}$, which indicates that $m \mid \gamma(q^{j}-1)$ and thus $\tau \mid j$. It follows that $\tau$ is the smallest positive integer satisfying \eqref{eq 3}, and
		$$c_{\ell m}(\gamma+\sigma m) = \{\gamma+\sigma m,(\gamma+\sigma m)q,\cdots,(\gamma+\sigma m)q^{\tau-1}\}.$$
		
		To prove the identity \eqref{eq 8}, first observe that the inclusion
		$$\bigcup_{\sigma=0}^{\ell-1}c_{\ell m}(\gamma+\sigma m) \subseteq \pi_{\ell m/m}^{-1}(c_{m}(\gamma))$$
		holds, as $\pi_{\ell m/m}(\gamma+\sigma m)=\gamma$ for any $0 \leq \sigma \leq \ell-1$. For $0 \leq \sigma_{1},\sigma_{2} \leq \ell-1$, if there is a positive integer $j$ such that
		$$(\gamma+\sigma_{1}m)q^{j} \equiv \gamma+\sigma_{2}m \pmod{\ell m},$$
		then $\gamma(q^{j}-1) \equiv m(\sigma_{1}q^{j}-\sigma_{2}) \pmod{\ell m}$, which implies that $m \mid \gamma(q^{j}-1)$, and $j$ is a multiple of $\tau$. Thus $\sigma_{1} = \sigma_{2}$. Consequently the cosets $c_{\ell m}(\gamma+\sigma m)$. $0 \leq \sigma \leq \ell-1$, are pairwise different. Now the identity \eqref{eq 8} just follows that 
		$$|\bigsqcup_{\sigma=0}^{\ell-1}c_{\ell m}(\gamma+\sigma m)| = \ell\cdot\tau =  |\pi_{\ell m/m}^{-1}(c_{m}(\gamma))|.$$
		
		\item[(2)] As $\gamma(q^{\tau}-1)$ is divisible by $m$ but not by $\ell m$, one obtains $v_{\ell}(\gamma) + v_{\ell}(q^{\tau}-1) = v_{\ell}(m)$. By Lemma \ref{lem 3} $\ell\tau$ is the smallest positive integer such that $\gamma q^{\ell\tau} \equiv \gamma \pmod{\ell m}$. Then we have
		$$c_{\ell m}(\gamma) = \{\gamma,\gamma q,\cdots,\gamma q^{\ell\tau-1}\}.$$
		Since $c_{\ell m}(\gamma) \subseteq \pi_{\ell m/m}^{-1}(c_{m}(\gamma))$, and
		$$|\pi_{\ell m/m}^{-1}(c_{m}(\gamma))| = \ell\cdot\tau =  |c_{\ell m}(\gamma)|,$$
		then
		$$\pi_{\ell m/m}^{-1}(c_{m}(\gamma)) = c_{\ell m}(\gamma) .$$
		It is obvious that in this case all the cosets $c_{\ell m}(\gamma+\sigma m)$, $0 \leq \sigma \leq \ell-1$, coincide.
	\end{itemize}
\end{proof}

On the other hand, if $\ell \nmid q^{\tau}-1$, then $q^{\tau}$ generates a nontrivial subgroup $\langle q^{\tau}\rangle$ in the multiplicative group $(\mathbb{Z}/\ell\mathbb{Z})^{\ast}$. We fix a set
$$\Sigma_{\tau} = \{\sigma_{1},\cdots,\sigma_{r}\} \subseteq \{1,\cdots,\ell-1\}$$
that represents the quotient group $(\mathbb{Z}/\ell\mathbb{Z})^{\ast}/ \langle q^{\tau}\rangle$, that is, $(\mathbb{Z}/\ell\mathbb{Z})^{\ast}/ \langle q^{\tau}\rangle$ can be written as
$$(\mathbb{Z}/\ell\mathbb{Z})^{\ast}/ \langle q^{\tau}\rangle = \{\sigma_{1}+\langle q^{\tau}\rangle,\cdots,\sigma_{r}+\langle q^{\tau}\rangle \}.$$
In particular, $r = [(\mathbb{Z}/\ell\mathbb{Z})^{\ast}:\langle q^{\tau}\rangle]$.

\begin{lemma}\label{lem 1}
	Let $m$ be a positive integer coprime to $q$, and $c_{m}(\gamma)$ be a $q$-cyclotomic coset modulo $m$ of size $\tau$. Assume that $\ell$ does not divide $q^{\tau}-1$. Choose a representative $\gamma_{0}$ of $\gamma+m\mathbb{Z} \in \mathbb{Z}/m\mathbb{Z}$ such that $v_{\ell}(\gamma_{0}) \geq v_{\ell}(m)+1$. Then the preimage of $c_{m}(\gamma)$ along the projection $\pi_{\ell m/m}$ is
	\begin{equation}\label{eq 5}
		\pi_{\ell m/m}^{-1}(c_{m}(\gamma)) = c_{\ell m}(\gamma_{0}) \sqcup (\bigsqcup_{\sigma \in \Sigma_{\tau}} c_{\ell m}(\gamma_{0}+\sigma m)).
	\end{equation}
	Moreover, the $q$-cyclotomic cosets modulo $\ell m$ on the RHS of \eqref{eq 5} are given by
	$$c_{\ell m}(\gamma_{0}) = \{\gamma_{0},\gamma_{0}q,\cdots,\gamma_{0}q^{\tau-1}\}$$
	and
	$$c_{\ell m}(\gamma_{0}+\sigma m)=\{\gamma_{0}+\sigma m,(\gamma_{0}+\sigma m)q,\cdots,(\gamma_{0}+\sigma m)q^{\mathrm{ord}_{\ell}(q^{\tau})\cdot\tau-1}\},$$
	for every $\sigma \in \Sigma_{\tau}$.
\end{lemma}

\begin{proof}
	First we explain that there is a representative $\gamma_{0}$ of $\gamma+m\mathbb{Z}$ such that $v_{\ell}(\gamma_{0}) \geq v_{\ell}(m)+1$. Since $\gamma q^{\tau} \equiv \gamma \pmod{m}$ and $\ell \nmid q^{\tau}-1$, then $v_{\ell}(\gamma) \geq v_{\ell}(m)$. If $v_{\ell}(\gamma) \geq v_{\ell}(m)+1$, setting $\gamma_{0} = \gamma$ meets the condition. Otherwise, we have $v_{\ell}(\gamma) = v_{\ell}(m)$, and thus there is a unique integer $1 \leq j \leq \ell-1$ such that
	$$jm \equiv - \gamma \pmod{\ell^{v_{\ell}(m)+1}},$$
	that is, $v_{\ell}(\gamma+jm) \geq v_{\ell}(m)+1$. Therefore $\gamma_{0} = \gamma + jm$ is what we seek for.
	
	Let $\gamma_{0}$ be a representative of $\gamma+m\mathbb{Z} \in \mathbb{Z}/m\mathbb{Z}$ such that $v_{\ell}(\gamma_{0}) \geq v_{\ell}(m)+1$. Then it is obvious that 
	$$c_{\ell m}(\gamma_{0}) = \{\gamma_{0},\gamma_{0}q,\cdots,\gamma_{0}q^{\tau-1}\}.$$
	For any $\sigma \in \Sigma_{\tau}$, since $v_{\ell}(\gamma_{0}+\sigma m) = v_{\ell}(m)$, then $\mathrm{ord}_{\ell}(q^{\tau})\cdot \tau$ is the smallest positive integer such that 
	$$(\gamma_{0}+\sigma m)q^{\mathrm{ord}_{\ell}(q^{\tau})\cdot\tau} \equiv \gamma_{0}+\sigma m \pmod{\ell m}.$$
	Therefore we have 
	$$c_{\ell m}(\gamma_{0}+\sigma m) = \{\gamma_{0}+\sigma m,(\gamma_{0}+\sigma m)q,\cdots,(\gamma_{0}+\sigma m)q^{\mathrm{ord}_{\ell}(q^{\tau})\cdot\tau-1}\}.$$
	
	Note that for any $\sigma_{1},\sigma_{2} \in \Sigma_{\tau}$, if there exists a positive integer $j$ such that 
	$$\gamma_{0}+ \sigma_{2}m \equiv (\gamma_{0} +\sigma_{1}m)q^{j} \pmod{\ell m},$$
	then
	$$\gamma_{0} \equiv \gamma_{0}q^{j} \pmod{m},$$
	which indicates that $\tau \mid j$ and $\sigma_{2}m \equiv \sigma_{1}mq^{j} \pmod{\ell m}$, or equivalently, $\sigma_{2} \equiv \sigma_{1}q^{j} \pmod{\ell}$. The assumption $\sigma_{1},\sigma_{2} \in \Sigma_{\tau}$ guarantees that $\sigma_{1} =\sigma_{2}$. Consequently, the cosets $c_{\ell m}(\gamma_{0})$ and $c_{\ell m}(\gamma_{0}+\sigma m)$, $\sigma \in \Sigma_{\tau}$, are pairwise different.
	
	Finally, as there are in total
	$$\tau + \sum_{\sigma \in \Sigma_{\tau}}\tau\cdot \mathrm{ord}_{\ell}(q^{\tau}) = \tau + \tau\cdot \mathrm{ord}_{\ell}(q^{\tau})\cdot \dfrac{\ell-1}{\mathrm{ord}_{\ell}(q^{\tau})} = \ell\tau$$
	elements in the disjoint union $c_{\ell m}(\gamma_{0}) \sqcup (\bigsqcup\limits_{\sigma \in \Sigma_{\tau}} c_{\ell m}(\gamma_{0}+\sigma m))$, then we obtain
	$$\pi_{\ell m/m}^{-1}(c_{m}(\gamma)) = c_{\ell m}(\gamma_{0}) \sqcup (\bigsqcup_{\sigma \in \Sigma_{\tau}} c_{\ell m}(\gamma_{0}+\sigma m)).$$
\end{proof}

\begin{remark}
	\begin{itemize}
		\item[(1)] In Lemma \ref{lem 4}, either $v_{\ell}(\gamma) + v_{\ell}(q^{\tau}-1) \geq v_{\ell}(m)+1$ or $v_{\ell}(\gamma) + v_{\ell}(q^{\tau}-1) = v_{\ell}(m)$ does not depend on the choice of the representative $\gamma$ of $\overline{\gamma}$, which can be seen as follow. If $\gamma$ is a representative of $\overline{\gamma}$ such that $v_{\ell}(\gamma) + v_{\ell}(q^{\tau}-1) = v_{\ell}(m)$, then $v_{\ell}(\gamma) = v_{\ell}(m) - v_{\ell}(q^{\tau}-1) < v_{\ell}(m)$ implies that for any $d \in \mathbb{Z}$
		$$v_{\ell}(\gamma+dm) = v_{\ell}(\gamma).$$
		\item [(2)] Also in Lemma \ref{lem 1} the disjoint union decomposition \eqref{eq 5} is independent of the choice of the representative $\gamma_{0}$. In fact, if $\gamma_{0}$ and $\gamma_{0}^{\prime}$ are two representatives of $\overline{\gamma}$ such that $v_{\ell}(\gamma_{0}), v_{\ell}(\gamma_{0}^{\prime}) \geq v_{\ell}(m)+1$, then $\ell m \mid \gamma_{0}-\gamma_{0}^{\prime}$ and therefore $\gamma_{0}$ and $\gamma_{0}^{\prime}$ represent the same element in $\mathbb{Z}/\ell m\mathbb{Z}$.
	\end{itemize}
\end{remark}

\begin{definition}
	Let $m$ be a positive integer coprime to $q$, and $c_{m}(\gamma)$ be a $q$-cyclotomic coset modulo $m$, with $| c_{m}(\gamma) | =\tau$.
	\begin{itemize}
		\item[(1)] If $\ell \nmid q^{\tau} -1$, the $q$-cyclotomic coset $c_{m}(\gamma)$ is said to be semi-splitting with respect to the extension $\mathbb{Z}/\ell m\mathbb{Z}: \mathbb{Z}/m\mathbb{Z}$.
		\item[(2)] Assume that $\ell \mid q^{\tau} -1$. If there is a representative (hence every representative) $\gamma$ of $\overline{\gamma}$ such that $v_{\ell}(\gamma) + v_{\ell}(q^{\tau}-1) \geq v_{\ell}(m)+1$, the $q$-cyclotomic coset $c_{m}(\gamma)$ is said to be splitting with respect to the extension $\mathbb{Z}/\ell m\mathbb{Z}: \mathbb{Z}/m\mathbb{Z}$; otherwise, $c_{m}(\gamma)$ is said to be nonsplitting or stable with respect to $\mathbb{Z}/\ell m\mathbb{Z}: \mathbb{Z}/m\mathbb{Z}$.
		\item[(3)] Let $n$ be a positive integer divisible by neither $p$ nor $\ell$, and $\mathfrak{c}=(c_{\ell^{i}n}(\gamma_{i}))_{i \in \mathbb{N}}$ be a sequence of $q$-cyclotomic cosets lying in the $\ell$-adic $q$-cyclotomic system $\mathcal{PC}_{n}^{\ell}$. For any $i_{0} \in \mathbb{N}$, we say that $\mathfrak{c}$ is semi-splitting (resp. splitting, nonsplitting) at degree $i_{0}$ if its $i_{0}$-th component $c_{\ell^{i_{0}}n}(\gamma_{i_{0}})$ is semi-splitting (resp. splitting, nonsplitting) with respect to $\mathbb{Z}/\ell^{i_{0}+1}n\mathbb{Z}: \mathbb{Z}/\ell^{i_{0}}n\mathbb{Z}$.
	\end{itemize}
\end{definition}

Let $\overline{\gamma} \in \mathbb{Z}/m\mathbb{Z}$. If there is a representative of $\overline{\gamma}$, say $\gamma$, is divisible by $\ell^{v_{\ell}(m)}$, then so is any representative $\gamma+\sigma m$, $\sigma \in \mathbb{Z}$, of $\overline{\gamma}$. Thus it makes sense to say whether an element in $\mathbb{Z}/m\mathbb{Z}$ is divisible by $\ell^{v_{\ell}(m)}$.

\begin{proposition}\label{prop 1}
	Let $m$ be a positive integer coprime to $q$, with $v_{\ell}(m)=v$.
	\begin{description}
		\item[(1)] If $\overline{\gamma} \in \mathbb{Z}/m\mathbb{Z}$ satisfies $\ell^{v} \mid \gamma$, then there is exactly one $q$-cyclotomic coset modulo $\ell m$ contained in $\widehat{\pi}_{\ell m/m}^{-1}(c_{m}(\gamma))$ whose elements are all divisible by $\ell^{v+1}$.
		\item[(2)] If $\overline{\gamma} \in \mathbb{Z}/m\mathbb{Z}$ is such that $\ell^{v} \nmid \gamma$, then any $q$-cyclotomic coset modulo $\ell m$ contained in $\widehat{\pi}_{\ell m/m}^{-1}(c_{m}(\gamma))$ consists of elements not divisible by $\ell^{v+1}$.
	\end{description}
\end{proposition}

\begin{proof}
	As $q$ is coprime to $\ell$, a $q$-cyclotomic coset modulo $\ell^{v+1}n$ consists of elements all divisible by $\ell^{v+1}$ if and only if it contains an element divisible by $\ell^{v+1}$. If $\ell \nmid q^{\tau}-1$ for $\tau = | c_{\ell^{v}n/q}(\gamma) |$, then by $v=v_{\ell}(m) < v_{\ell}(\gamma q^{\tau}-\gamma) = v_{\ell}(\gamma)$ one has $\ell^{v} \mid \gamma$. According to Lemma \ref{lem 1}, $\widehat{\pi}_{\ell m/m}^{-1}(c_{m}(\gamma))$ can be decomposed as
	$$\pi_{\ell m/m}^{-1}(c_{m}(\gamma)) = c_{\ell m}(\gamma_{0}) \sqcup (\bigsqcup_{\sigma \in \Sigma_{\tau}} c_{\ell m}(\gamma_{0}+\sigma m)),$$
	where $\gamma_{0}$ is a representative of $\gamma$ which is divisible by $\ell^{v+1}$. Notice that for any $\sigma \in \Sigma_{\tau} \subseteq \{1,\cdots,\ell-1\}$, $v_{\ell}(\gamma_{0}+\sigma m)=v_{\ell}(m)=v$. Hence $c_{\ell m}(\gamma_{0})$ is the unique coset contained in $\widehat{\pi}_{\ell m/m}^{-1}(c_{m}(\gamma))$ whose elements are divisible by $\ell^{v+1}$.
	
	Next we suppose that $\ell \mid q^{\tau}-1$. By Lemma \ref{lem 4} we have 
	\begin{equation*}
		\pi_{\ell m/m}^{\-1}(c_{m}(\gamma)) = \left\{
		\begin{array}{lcl}
			\bigsqcup\limits_{\sigma=0}^{\ell-1}c_{\ell m}(\gamma+\sigma m), \ v_{\ell}(\gamma)+v_{\ell}(q^{\tau}-1) \geq v+1;\\
			c_{\ell m}(\gamma), \ v_{\ell}(\gamma)+v_{\ell}(q^{\tau}-1) = v.
		\end{array} \right.
	\end{equation*}
	If $\ell^{v} \mid \gamma$, then $v_{\ell}(\gamma)+v_{\ell}(q^{\tau}-1) \geq v+1$. There is exactly one $\sigma_{0} \in \{0,\cdots,\ell-1\}$ such that
	$$\sigma_{0}m \equiv -\gamma \pmod{\ell^{v+1}},$$
	therefore among $c_{\ell m}(\gamma+\sigma m)$, $0 \leq \sigma \leq \ell-1$, there is exactly one coset whose elements are divisible by $\ell^{v+1}$. Otherwise, if $\ell^{v} \nmid \gamma$, then for any $\sigma \in \mathbb{Z}$
	$$v_{\ell}(\gamma+\sigma m) = v_{\ell}(\gamma) < v.$$
	Thus in either case, the cosets lying in $\widehat{\pi}_{\ell m/m}^{-1}(c_{m}(\gamma))$ consist of elements that are not divisible by $\ell^{v+1}$.
\end{proof}

For any $q$-cyclotomic coset $c_{n}(\gamma)$ modulo $n$, as $\ell \nmid n$ then $\ell^{v_{\ell}(n)} \mid \gamma$ trivially holds. Applying Proposition \ref{prop 1} inductively yields that there is exactly one sequence $(c_{\ell^{i}n/q}(\gamma_{i}))_{i \in \mathbb{N}}$ in $\mathcal{PC}_{n}^{\ell}(\gamma)$ which satisfies that $\gamma_{i}$ is divisible by $\ell^{i+1}$ for every $i \in \mathbb{N}$, and consequently is splitting at every degree. (A more detailed proof is included in the proofs of Theorem \ref{thm 3} and \ref{thm 4}.) This fact motivates the notions of cyclotomic $\ell$-adic class and cyclotomic $\ell$-adic integer, which will play the central roles in the classification of the elements in $\mathcal{PC}_{n}^{\ell}(\gamma)$ and the characterization of its structure presented in the next two sections.

Let
$$\mathbb{Z}_{\ell} = \varprojlim\limits_{i \in \mathbb{N}} \mathbb{Z}/\ell^{i}\mathbb{Z}$$ 
be the ring of $\ell$-adic integers. Define a homomorphism of abelian groups 
$$\varphi_{n}^{\ell}: \mathbb{Z} \rightarrow \mathbb{Z}_{\ell}; \ \gamma \mapsto -\frac{\gamma}{n}.$$ 
Since $\ell$ does not divide $n$, $-\frac{\gamma}{n}$ lies in $\mathbb{Z}_{\ell}$ and the homomorphism $\varphi_{n}^{\ell}$ is well-defined. Further, consider the composition
$$\overline{\varphi}_{n}^{\ell}: \mathbb{Z} \rightarrow \mathbb{Z}_{\ell} \rightarrow \mathbb{Z}_{\ell}/\mathbb{Z}; \ \gamma \mapsto \overline{-\dfrac{\gamma}{n}} = -\dfrac{\gamma}{n} + \mathbb{Z}.$$
It is obvious that the kernel of $\overline{\varphi}_{n}^{\ell}$ is $n\mathbb{Z}$, hence we obtain an induced homomorphism, which is still denoted by $\overline{\varphi}_{n}^{\ell}$:
$$\overline{\varphi}_{n}^{\ell}: \mathbb{Z}/n\mathbb{Z} \rightarrow \mathbb{Z}_{\ell}/\mathbb{Z}; \ \overline{\gamma} \mapsto \overline{-\dfrac{\gamma}{n}}.$$

\begin{definition}
	Let $\overline{\gamma} \in \mathbb{Z}/n\mathbb{Z}$. Define the cyclotomic $\ell$-adic class associated to $\overline{\gamma}$ to be the image $\overline{\varphi}_{n}^{\ell}(\overline{\gamma})$. Each $\ell$-adic integer lying in the class $\overline{\varphi}_{n}^{\ell}(\overline{\gamma})$ is called a cyclotomic $\ell$-adic integer over $\overline{\gamma}$.
\end{definition}

As each $\ell$-adic integer can be identified with its $\ell$-adic expansion, the image $\varphi_{n}^{\ell}(\gamma)$ of any integer $\gamma$ can be presented uniquely as a power series
\begin{equation}\label{eq 2}
	\varphi_{n}^{\ell}(\gamma) = \sum_{k=0}^{\infty}\phi_{k}\cdot\ell^{k},
\end{equation} 
where $\phi_{k} \in [0,\ell-1]$. Furthermore, the quotient group $\mathbb{Z}_{\ell}/\mathbb{Z}$ can be identified with the quotient of the abelian group of $\ell$-adic power series modulo the subgroup of the series corresponding to rational integers. Under such identifications, $\varphi_{n}^{\ell}(\gamma)$ and $\overline{\varphi}_{n}^{\ell}(\overline{\gamma})$ can be equivalently characterized by the following property.

\begin{proposition}\label{prop 2}
	\begin{itemize}
		\item[(1)] For any integer $\gamma$, the image $\varphi_{n}^{\ell}(\gamma) = \sum\limits_{k=0}^{\infty}\phi_{k}\cdot\ell^{k}$ is the unique $\ell$-adic power series satisfying that for all $i \in \mathbb{N}$
		$$\ell^{i+1} \mid \gamma+ n\cdot \varphi_{n}^{\ell}(\gamma)_{\leq i},$$
		where $\varphi_{n}^{\ell}(\gamma)_{\leq i} = \sum\limits_{k=0}^{i}\phi_{k}\cdot\ell^{k}$ is the finite sum of the first $i+1$ terms of $\varphi_{n}^{\ell}(\gamma)$.
		\item[(2)] The map $\varphi_{n}^{\ell}$ induces a bijection from the set of representatives of $\overline{\gamma}$ in $\mathbb{Z}$ onto the set of cyclotomic $\ell$-adic integers over $\overline{\gamma}$.
	\end{itemize}
\end{proposition}

\begin{proof}
	Since $n$ is coprime to $\ell$, by induction it can be checked that there exists a unique series $\sum\limits_{k=0}^{\infty}a_{k}\cdot\ell^{k}$ such that 
	\begin{equation}\label{eq 4}
		\ell^{i+1} \mid \gamma + n\cdot \sum_{k=0}^{i}a_{k}\cdot\ell^{k}, \ \forall i \in \mathbb{N}.
	\end{equation}
	Notice that \eqref{eq 4} is equivalent to $\gamma + n\cdot \sum\limits_{k=0}^{\infty}a_{k}\cdot\ell^{k} = 0$ in the ring $\mathbb{Z}_{\ell}$, therefore we have 
	$$\sum\limits_{k=0}^{\infty}a_{k}\cdot\ell^{k} = -\dfrac{\gamma}{n} = \varphi_{n}^{\ell}(\gamma).$$
	The conclusion $(2)$ follows from $(1)$ directly.
\end{proof}

Following from $(2)$ of Proposition \ref{prop 2}, for a representative $\gamma$ of $\overline{\gamma} \in \mathbb{Z}/n\mathbb{Z}$, the image $\varphi_{n}^{\ell}(\gamma)$ is called the cyclotomic $\ell$-adic integer associated to $\gamma$.

\begin{remark}
	From the construction of the homomorphisms $\varphi_{n}^{\ell}$ and $\overline{\varphi}_{n}^{\ell}$, it can be seen that they are independent of $q$.
\end{remark}

\section{$\ell$-adic $q$-cyclotomic system: odd prime $\ell$}\label{sec 3}
The present and the next sections are devoted to classifying the elements in the $\ell$-adic $q$-cyclotomic system $\mathcal{PC}_{n}^{\ell}$ and characterizing its structure, in the cases where $\ell$ is an odd prime and where $\ell=2$ respectively. First in this section we assume that $\ell$ is an odd prime different from $p$, and $n$ is a positive integer that is not divisible by $p$ and by $\ell$.

By Lemma \ref{lem 6} there is a disjoint union decomposition
$$\mathcal{PC}_{n}^{\ell} = \bigsqcup_{\gamma \in \mathcal{CR}_{n}}\mathcal{PC}_{n}^{\ell}(\gamma),$$
where $\mathcal{CR}_{n}$ is any full set of representatives of $\mathcal{C}_{n}$. Thus we only need to understand the $\ell$-adic $q$-cyclotomic system $\mathcal{PC}_{n}^{\ell}(\gamma)$ over $c_{n}(\gamma)$ for every $\gamma \in \mathcal{CR}_{n}$.

Let $c_{n}(\gamma)$ be a $q$-cyclotomic coset modulo $n$, with size $|c_{n}(\gamma)|=\tau$. Assume that $\ell$ divides $q^{\tau}-1$. Let $\varphi_{n}^{\ell}(\gamma) = \sum\limits_{k=0}^{\infty}\phi_{k}\cdot\ell^{k}$ be the cyclotomic $\ell$-adic integer associated to $\gamma$. In this section, as $\ell$ and $n$ are fixed, we will simply write $\varphi(\gamma)$ for $\varphi_{n}^{\ell}(\gamma)$. For any $m \in \mathbb{N}$ the coefficient $\phi_{k}$ of $\varphi_{n}^{\ell}(\gamma)$ lies in the set $[0,\ell-1]$. We define $\Sigma_{m} = [0,\ell-1]\setminus\{\phi_{m}\}$, and write it as
$$\Sigma_{m} = \{\sigma_{1}^{m},\cdots,\sigma_{\ell-1}^{m}\}$$
with $\sigma_{m,1} < \cdots < \sigma_{m,\ell-1}$. Then for any pair $(m,a)$ consisting of $m \in \mathbb{N}$ and $a \in [1,\ell-1]$ we define a rational integer $U_{n}^{\ell}(\gamma;m,a)$ by the sequence $U_{n}^{\ell}(\gamma;m,a) = \sum\limits_{k=0}^{m}u_{k}(m,a)\cdot\ell^{k}$,
\begin{equation*}
	u_{k}(m,a) = 
	\left\{
	\begin{array}{lcl}
		\phi_{k}, \ 0 \leq k \leq m-1;\\
		\sigma_{a}^{m}, \ k=m.
	\end{array} \right.
\end{equation*}
If the parameters $n$ and $\ell$ are clear from the context, we will omit the superscript and write $U_{n}^{\ell}(\gamma;m,a) = U(\gamma;m,a)$.

\begin{theorem}\label{thm 3}
	Let $c_{n}(\gamma)$ be a $q$-cyclotomic coset modulo $n$ with size $\tau$. Assume that $\ell$ divides $q^{\tau}-1$. 
	\begin{description}
		\item[(1)] The cyclotomic $\ell$-adic integer $\varphi(\gamma)$ gives rise to an element
		$$\mathfrak{c}_{\infty}(\gamma) = (c_{\ell^{i}n}(\gamma+n\cdot\varphi(\gamma)_{\leq i-1}))_{i \in \mathbb{N}} \in \mathcal{PC}_{n}^{\ell}(\gamma).$$
		Here $\varphi(\gamma)_{\leq 1}$ is set to be $0$. The sequence $\mathfrak{c}_{\infty}(\gamma)$ is the unique splitting element in $\mathcal{PC}_{n}(\gamma)$.
		\item[(2)] For any $m \in \mathbb{N}$ and $a \in [1,\ell-1]$, the integer $U(\gamma;m,a)$ gives rise to $\ell^{v-1}$ elements
		$$\mathfrak{c}_{m,a,t}(\gamma) = (c_{\ell^{i}n}(\gamma+n(U(\gamma;m,a)+\ell^{m+1}t)))_{i \in \mathbb{N}} \in \mathcal{PC}_{n}^{\ell}(\gamma),$$
		where $v = v_{\ell}(q^{\tau}-1)$ and $t \in [0,\ell^{v-1}-1]$.
	\end{description}
	Furthermore, all the elements in $\mathcal{PC}_{n}^{\ell}(\gamma)$ are obtained by the above construction unrepeatedly.
\end{theorem}

\begin{proof}
	It is clear that the sequences given in $(1)$ and $(2)$ lie in $\mathcal{PC}_{n}(\gamma)$. By Proposition \ref{prop 2} one has
	$$v_{\ell}(\gamma+n\cdot\varphi(\gamma)_{\leq i-1}) \geq i$$
	for every $i \in \mathbb{N}$, therefore $\mathfrak{c}_{\infty}(\gamma)$ is splitting. On the other hand, by the definition of $U(\gamma;m,a)$
	$$v_{\ell}(\gamma+n(U(\gamma;m,a)+\ell^{m+1}t))=m$$
	for any $t \in [0,\ell^{v-1}-1]$, which implies that $\mathfrak{c}_{m,a,t}(\gamma)$ is splitting at degree $i \leq m+v-1$ while is stable at degree $i \geq m+v$. Consequently the sequences in $(1)$ and $(2)$ are pairwise distinct.
	
	It remains to prove that the sequences in $(1)$ and $(2)$ give all the elements in $\mathcal{PC}_{n}^{\ell}(\gamma)$. Let $\mathfrak{c}=(c_{\ell^{i}n}(\gamma_{i}))_{i \in \mathbb{N}} \in \mathcal{PC}_{n}^{\ell}(\gamma)$. Applying Lemma \ref{lem 4} successively yields that the representatives $\gamma_{i}$ can be chosen to be in the form
	$$\gamma_{i} = \gamma + n\cdot Z_{\leq i-1},$$
	where $Z_{\leq i-1}$ is the finite sum of the first $i$ terms of a power series $Z = \sum\limits_{k=0}^{\infty}z_{k}\cdot\ell^{k}$ for $i \geq 1$ and $Z_{\leq -1}=0$. If $Z = \varphi(\gamma)$ then $\mathfrak{c}$ is exactly $c_{\infty}(\gamma)$. Otherwise there exists a minimal nonnegative integer $m$ such that $z_{m} \neq \phi_{m}$, then the component at degree $m$ of $\mathfrak{c}$ is $c_{\ell^{m}n}(\gamma+n\cdot\varphi_{n}^{\ell}(\gamma)_{\leq m-1})$. By Lemma \ref{lem 4} one obtains
	$$\pi_{\ell^{m+1}n/\ell^{}n}^{-1}(c_{\ell^{m}n}(\gamma+n\cdot\varphi_{n}^{\ell}(\gamma)_{\leq m-1})) = \bigsqcup_{a=0}^{\ell-1}c_{\ell^{m+1}n}(\gamma+n\cdot(\varphi_{n}^{\ell}(\gamma)_{\leq m-1}+a\cdot\ell^{m})).$$
	As $z_{m} \neq \phi_{m}$, $z_{m} \in [0,\ell-1] \setminus \{\phi_{m}\} = \Sigma_{m}$, that is, there is an integer $1 \leq a \leq \ell-1$ such that 
	\begin{equation}\label{eq 9}
		Z_{\leq m} = U(\gamma;m,a).
	\end{equation}
	By Proposition \ref{prop 1} for any power series $Z = \sum\limits_{k=0}^{\infty}z_{k}\cdot\ell^{k}$ satisfying \eqref{eq 9} it holds that
	$$v_{\ell}(\gamma+n\cdot Z_{\leq i}) = m$$
	for all $i \geq m$. Consequently the sequence $\mathfrak{c}$ is splitting at degree $i \leq m+v-1$, while is stable at degree $i \geq m+v$. It follows that   
	$$Z = U(\gamma;m,a) + \sum_{j=1}^{v-1}t_{j}\ell^{m+j} = U(\gamma;m,a)+\ell^{m+1}\cdot\sum_{j=0}^{v-2}t_{j}\ell^{j}$$
	for some $(v-1)$-tuple $(t_{1},\cdots,t_{v-1}) \in [0,\ell-1]^{v-1}$. Set $t = \sum\limits_{j=0}^{v-2}t_{j}\ell^{j}$. Then $t \in [0,\ell^{v-1}-1]$ and
	$$Z = U(\gamma;m,a)+\ell^{m+1}t,$$
	which yields $\mathfrak{c} = c_{m,a,t}(\gamma)$.
\end{proof}

Further, the size of every component of the sequences in $\mathcal{PC}_{n}^{\ell}(\gamma)$ are computed in the corollary below.

\begin{corollary}\label{coro 1}
	Let the notations be given as in Theorem \ref{thm 3}. Then
	\begin{description}
		\item[(1a)] Every component of the sequence $\mathfrak{c}_{\infty}(\gamma)$ has the size $\tau = | c_{n}(\gamma) |$.
		\item[(1b)] For any pair $(m,a) \in \mathbb{N} \times [1,\ell-1]$, the size of each component of $\mathfrak{c}_{m,a,t}(\gamma)$ is given by
		\begin{equation*}
			| c_{\ell^{i}n}(\gamma+n(U(\gamma;m,a)+\ell^{m+1}t)) | =\left\{
			\begin{array}{lcl}
				\tau, \quad 0 \leq i \leq m+v;\\
				\ell^{i-m-v}\cdot\tau, \quad i \geq m+v+1.
			\end{array} \right.
		\end{equation*}
	\end{description}
\end{corollary}

\begin{proof}
	The sequence $\mathfrak{c}_{\infty}(\gamma)$ is splitting, and by the proof of Theorem \ref{thm 3} for any $(m,a) \in \mathbb{N}\times[1,\ell-1]$ the sequence $\mathfrak{c}_{m,a,t}(\gamma)$ is splitting at degree $i \leq m+v-1$, while is stable at degree $i \geq m+v$. Then the conclusions follow from Lemma \ref{lem 4} directly.
\end{proof}

Next we turn to the case that $\ell$ does not divide $q^{\tau}-1$. In this case $q^{\tau}$ generates a nontrivial cyclic subgroup $\langle q^{\tau}\rangle$ in $(\mathbb{Z}/\ell\mathbb{Z})^{\ast}$ with index $r = \frac{\ell-1}{\mathrm{ord}_{\ell}(q^{\tau})}$. Recall that the set 
$$\Sigma_{\tau} = \{\sigma_{1},\cdots,\sigma_{r}\} \subseteq [1,\ell-1]$$
is defined to represents the quotient group $(\mathbb{Z}/\ell\mathbb{Z})^{\ast}/\langle q^{\tau}\rangle$. We further require that $\sigma_{1} < \cdots < \sigma_{r}$. To any pair $(m,a)$ consisting of a nonnegative integer $m$ and an integer $a \in [1,r]$ we attach a rational integer $U_{n}^{\ell}(\gamma;m,a) = \sum\limits_{k=0}^{m}u_{k}(m,a)\cdot\ell^{k}$ by
\begin{equation*}
    u_{k}(m,a) =\left\{
	\begin{array}{lcl}
		\phi_{k}, \quad 0 \leq k \leq m-1;\\
		(\phi_{m}+\sigma_{a}) \pmod{\ell}, \quad k = m.
	\end{array} \right.
\end{equation*}
Here $(\phi_{m}+\sigma_{a}) \pmod{\ell}$ stands for the unique integer $\sigma$ such that $0 \leq \sigma \leq \ell-1$ and 
$$\sigma \equiv \phi_{m}+\sigma_{a} \pmod{\ell}.$$
Also, when the parameters $n$ and $\ell$ are clear, we denote $U_{n}^{\ell}(\gamma;m,a)$ simply by $U(\gamma;m,a)$.

The elements of $\mathcal{PC}_{n}^{\ell}(\gamma)$ and the sizes of their components, in the case where $\ell \nmid q^{\tau}-1$, are given as follows. The proofs are similar to that of Theorem \ref{thm 3} and Corollary \ref{coro 1}.

\begin{theorem}\label{thm 4}
	Let $c_{n}(\gamma)$ be a $q$-cyclotomic coset modulo $n$ with size $\tau$. Assume that $\ell$ does not divide $q^{\tau}-1$. 
	\begin{description}
		\item[(1)] The cyclotomic $\ell$-adic integer $\varphi(\gamma)$ gives rise to an element
		$$\mathfrak{c}_{\infty}(\gamma) = (c_{\ell^{i}n}(\gamma + n\cdot \varphi(\gamma)_{\leq i-1}))_{i \in \mathbb{N}} \in \mathcal{PC}_{n}^{\ell}(\gamma).$$
		Here $\varphi(\gamma)_{\leq 1}$ is set to be $0$. The sequence $\mathfrak{c}_{\infty}(\gamma)$ is the unique splitting element in $\mathcal{PC}_{n}(\gamma)$.
		\item[(2)] For each $m \in \mathbb{N}$ and $a \in [1,r]$, the integer $U(\gamma;m,a)$ gives rise to $\ell^{v-1}$ elements
		$$\mathfrak{c}_{m,a,t}(\gamma) = (c_{\ell^{i}n}(\gamma + n(U(\gamma;m,a)+\ell^{m+1}t)))_{i \in \mathbb{N}} \in \mathcal{PC}_{n}^{\ell}(\gamma),$$
		where $v = v_{\ell}(q^{\tau\cdot\mathrm{ord}_{\ell}(q^{\tau})}-1)$ and $t \in [0,\ell^{v-1}-1]$.
	\end{description}
	Furthermore, all elements in $\mathcal{PC}_{n}^{\ell}(\gamma)$ are obtained by the above construction unrepeatedly.
\end{theorem}

\begin{corollary}\label{coro 3}
	Let the notations be given as in Theorem \ref{thm 4}. Then
	\begin{description}
		\item[(1)] Every component of the sequence $\mathfrak{c}_{\infty}(\gamma)$ has the size $\tau = | c_{n}(\gamma) |$.
		\item[(2)] For any pair $(m,a) \in \mathbb{N} \times [1,r]$, the size of each component of $\mathfrak{c}_{m,a,t}(\gamma)$ is given by
		\begin{equation*}
			| c_{\ell^{i}n}(\gamma + n(U(\gamma;m,a)+\ell^{m+1}t)) | =\left\{
			\begin{array}{lcl}
				\tau, \quad 0 \leq i \leq m;\\
				\mathrm{ord}_{\ell}(q^{\tau})\cdot\tau, \quad m+1 \leq i \leq m+v;\\
				\ell^{i-m-v}\cdot\mathrm{ord}_{\ell}(q^{\tau})\cdot\tau, \quad i \geq m+v+1.
			\end{array} \right.
		\end{equation*}
	\end{description}
\end{corollary} 

Following from Theorem \ref{thm 3} and \ref{thm 4}, there is a unique element in $\mathcal{PC}_{n}^{\ell}(\gamma)$, the sequence $\mathfrak{c}_{\infty}(\gamma)$ induced by the cyclotomic $\ell$-adic integer $\varphi(\gamma)$, that is splitting or semi-splitting at every degree. We call $\mathfrak{c}_{\infty}(\gamma)$ the principal sequence over $c_{n}(\gamma)$, and denote
$$\mathcal{PC}_{n}^{\ell}(\gamma;\infty) = \{\mathfrak{c}_{\infty}(\gamma)\}.$$

The other elements in $\mathcal{PC}_{n}^{\ell}(\gamma)$ are called the stable sequences over $c_{n}(\gamma)$. For each $m \in \mathbb{N}$ we set
\begin{equation*}
	v =\left\{
	\begin{array}{lcl}
		v_{\ell}(q^{\tau}-1), \quad \mathrm{if} \ \ell \mid q^{\tau}-1;\\
		v_{\ell}(q^{\tau\cdot\mathrm{ord}_{\ell}(q^{\tau})}-1), \quad \mathrm{if} \ \ell \nmid q^{\tau}-1,
	\end{array} \right.
\end{equation*}
and define 
\begin{equation*}
	\mathcal{PC}_{n}^{\ell}(\gamma;m) =\left\{
	\begin{array}{lcl}
		\{\mathfrak{c}_{m,a,t}(\gamma) \ | \ a \in [1,\ell-1], \ t \in [0,\ell^{v-1}-1]\}, \quad \mathrm{if} \ \ell \mid q^{\tau}-1;\\
		\{\mathfrak{c}_{m,a,t}(\gamma) \ | \ a \in [1,r], \ t \in [0,\ell^{v-1}-1]\}, \quad \mathrm{if} \ \ell \nmid q^{\tau}-1.
	\end{array} \right.
\end{equation*}
The elements in $\mathcal{PC}_{n}^{\ell}(\gamma;m)$ are said to be of quasi-stable degree $m$ and of stable degree $m+v$. Then Theorem \ref{thm 3} and \ref{thm 4} can be rephrased uniformly as follow.

\begin{corollary}
	The $\ell$-adic $q$-cyclotomic system $\mathcal{PC}_{n}^{\ell}(\gamma)$ over $c_{n}(\gamma)$ can be naturally classified via the quasi-stable degree as
	$$\mathcal{PC}_{n}^{\ell}(\gamma) = \bigsqcup_{m \in \mathbb{N}\cup\{\infty\}}\mathcal{PC}_{n}^{\ell}(\gamma;m).$$
\end{corollary}

Summarizing the results obtained above, we are ready to state the theorem characterizing the structure of the $\ell$-adic $q$-cyclotomic system $\mathcal{PC}_{n}^{\ell}$.

\begin{theorem}\label{thm 7}
	Let $q$ be a power of a prime $p$, $\ell$ be an odd prime different than $p$, and $n$ be a positive integer not divisible by $p$ and by $\ell$. Then the $\ell$-adic $q$-cyclotomic system $\mathcal{PC}_{n}^{\ell}$ can be decomposed as
	$$\mathcal{PC}_{n}^{\ell} = \bigsqcup_{\gamma \in \mathcal{CR}_{n}}\left(\bigsqcup_{m \in \mathbb{N}\cup\{\infty\}}\mathcal{PC}_{n}^{\ell}(\gamma;m)\right).$$
	Furthermore, we have
	\begin{description}
		\item[(1)] The sequences lying in $\mathcal{PC}_{n}^{\ell}(\gamma)_{\infty}$, $\gamma \in \mathcal{CR}_{n}$, are exactly the sequences that are splitting or semi-splitting.
		\item[(2)] Let $\gamma \in \mathcal{CR}_{n}$ with $\tau = c_{n}(\gamma)$. If $\ell \mid q^{\tau}-1$, then for any $m \in \mathbb{N}$, each sequence lying in $\mathcal{PC}_{n}^{\ell}(\gamma;m)$ is splitting at degree $i \leq m+v-1$, while is stable at degree $i \geq m+v$. If $\ell \nmid q^{\tau}-1$, for any $m \in \mathbb{N}$, each sequence in $\mathcal{PC}_{n}^{\ell}(\gamma;m)$ is semi-splitting at degree $i \leq m-1$, splitting at degree $m \leq i \leq m+v-1$, and non-splitting at degree $i \geq m+v$.
	\end{description}
\end{theorem}

Finally we conclude this section with a concrete example.

\begin{examples}
	Let $q$ be a power of a prime $p$, and $\ell$ be an odd prime different than $p$. As there is only one $q$-cyclotomic coset modulo $1$, that is, $c_{1}(0) = \{0\}$, therefore $\mathcal{PC}_{1}^{\ell} = \mathcal{PC}_{1}^{\ell}(0)$. The cyclotomic $\ell$-adic integer associated to $0$ is just $0$.
	
	If $\ell \mid q-1$, for any $m \in \mathbb{N}$
	$$\Sigma_{m} = \{0,1,\cdots,\ell-1\} \setminus \{0\} = \{1,\cdots,\ell-1\},$$
	which implies that for $a \in \{1,\cdots,\ell-1\}$
	$$U(0;m,a) = a\ell^{m}.$$
	It follows from Theorem \ref{thm 3} that $\varphi(0) = 0$ induces the sequence
	$$\mathfrak{c}_{\infty}(0) = (c_{\ell^{i}}(0))_{i \in \mathbb{N}} \in \mathcal{PC}_{1}^{\ell},$$
	and for any $(m,a) \in \mathbb{N}\times [1,\ell-1]$, $U(0;m,a)$ induces $\ell^{v-1}$ sequences
	$$\mathfrak{c}_{m,a,t}(0) = (c_{\ell^{i}}(\ell^{m}(a+lt)))_{i \in \mathbb{N}} \in \mathcal{PC}_{1}^{\ell},$$
	where $v = v_{\ell}(q-1)$ and $t \in [0,\ell^{v-1}-1]$. Such sequences vary all the elements in $\mathcal{PC}_{1}^{\ell}$.
	
	Moreover, every component $c_{\ell^{i}/q}(0)$ of $\mathfrak{c}_{\infty}(0)$ has size $1$, and the size of each component of $\mathfrak{c}_{m,a,t}(0)$ is given by
	\begin{equation*}
		|c_{\ell^{i}}(\ell^{m}(a+lt))| =
		\left\{
		\begin{array}{lcl}
			1, \ \mathrm{if} \quad 0 \leq i \leq m+v;\\
			\ell^{i-m-v}, \quad \mathrm{if} \ i\geq m+v+1.
		\end{array} \right.
	\end{equation*}
	
	If $\ell \nmid q-1$, we fix a primitive root $\xi$ modulo $\ell$, and set $\xi_{a} = \xi^{a} \pmod{\ell}$ for $1 \leq a \leq r = \frac{\ell-1}{\mathrm{ord}_{\ell}(q)}$. Then the set $\Sigma_{1}$ can be chosen to be
	$$\Sigma_{1} = \{\xi_{a}\ |\ 0 \leq a \leq r-1\}.$$
	For any $(m,a) \in \mathbb{N}\times [1,r]$ 
	$$U(0;m,a) = \xi_{a}\ell^{m}.$$
	It follows from Theorem \ref{thm 3} that $\varphi(0) = 0$ induces the sequence
	$$\mathfrak{c}_{\infty}(0) = (c_{\ell^{i}}(0))_{i \in \mathbb{N}} \in \mathcal{PC}_{1}^{\ell},$$
	and for any $(m,a) \in \mathbb{N}\times [1,r]$, $U(0;m,a)$ induces $\ell^{v-1}$ sequences
	$$\mathfrak{c}_{m,a,t}(0) = (c_{\ell^{i}}(\ell^{m}(\xi_{a}+lt)))_{i \in \mathbb{N}} \in \mathcal{PC}_{1}^{\ell},$$
	where $v = v_{\ell}(q^{\mathrm{ord}_{\ell}(q)}-1)$ and $t \in [0,\ell^{v-1}-1]$. Such sequences vary all the elements in $\mathcal{PC}_{1}^{\ell}$.
	
	Moreover, every component $c_{\ell^{i}}(0)$ of $\mathfrak{c}_{\infty}(0)$ has size $1$, and the size of each component of $\mathfrak{c}_{m,a,t}(0)$ is given by
	\begin{equation*}
		|c_{\ell^{i}}(\ell^{m}(\xi_{a}+lt))| =
		\left\{
		\begin{array}{lcl}
			1, \ \mathrm{if} \ 0 \leq i \leq m;\\
			\mathrm{ord}_{\ell}(q), \ \mathrm{if} \ m+1 \leq i \leq m+v;\\
			\ell^{N-m-v}\cdot\mathrm{ord}_{\ell}(q), \ \mathrm{if} \ i \geq m+v+1.
		\end{array} \right.
	\end{equation*}
\end{examples}

\section{$\ell$-adic $q$-cyclotomic system: $\ell = 2$}\label{sec 4}
In this section we deal with the case that $\ell=2$. The main results have been contained in \cite{Zhu}. To bring both the cases where $\ell$ is an odd prime number and where $\ell=2$ under a unified framework, we will state these results without proving them. We refer the readers to Section $3.3.2$ in \cite{Zhu} for detailed proofs. Throughout this section it is assumed that $p$ is an odd prime, $q = p^{e}$ where $e \in \mathbb{N}^{+}$, and $n$ is a positive odd integer not divisible by $p$.

Let $c_{n}(\gamma)$ be a $q$-cyclotomic coset modulo $n$, with size $\tau$, and let 
$$\varphi(\gamma) = \varphi_{n}^{2}(\gamma) = \sum_{k=0}^{\infty}\phi_{k}\cdot 2^{k}$$
be the cyclotomic $2$-adic integer associated to $\gamma$. For any $m \in \mathbb{N}$ set
$$\Sigma_{m} = \{0,1\} \setminus \{\phi_{m}\} = \{1-\phi_{m}\},$$
and define a rational integer $U_{n}^{\ell}(\gamma;m)$ by the sequence $U_{n}^{\ell}(\gamma;m)= \sum\limits_{k=0}^{m}u_{k}(m)\cdot 2^{k}$, where 
\begin{equation*}
	u_{k}(m)=
	\left\{
	\begin{array}{lcl}
		\phi_{k}, \ 0 \leq k \leq m-1;\\
		1-\phi_{m}, \ k = m.
	\end{array} \right.
\end{equation*}
Still, we write $U_{n}^{\ell}(\gamma;m)$ simply as $U(\gamma;m)$ in this section as $n$ and $\ell=2$ are clear from the context.

The elements in the cyclotomic system $\mathcal{PC}_{n}^{2}(\gamma)$ and the sizes of their components are determined by the theorems and corollary below. 

\begin{theorem}\label{thm 5}
	Let $c_{n}(\gamma)$ be a $q$-cyclotomic coset modulo $n$ with $| c_{n}(\gamma) | =\tau$. Assume that $q^{\tau} \equiv 1 \pmod{4}$. Then
	\begin{description}
		\item[(1)] The cyclotomic $2$-adic integer $\varphi(\gamma)$ gives rise to an element
		$$\mathfrak{c}_{\infty}(\gamma) = (c_{2^{i}n}(\gamma + n\cdot \varphi(\gamma)_{\leq i-1}))_{i \in \mathbb{N}} \in \mathcal{PC}_{n}^{2}(\gamma).$$
		Here $\varphi(\gamma)_{\leq -1}$ is set to be $0$. The sequence $\mathfrak{c}_{\infty}(\gamma)$ is the unique splitting element in $\mathcal{PC}_{n}^{2}(\gamma)$.
		\item[(2)] For each $m \in \mathbb{N}$, the integer $U(\gamma;m)$ gives rise to $2^{v-1}$ elements
		$$\mathfrak{c}_{m,t}(\gamma) = (c_{2^{i}n}(\gamma + n(U(\gamma;m)+2^{m+1}t)))_{i \in \mathbb{N}} \in \mathcal{PC}_{n}^{2}(\gamma),$$
		where $v = v_{2}(q^{\tau}-1)$ and $t \in [0,2^{v-1}-1]$.
	\end{description}
	Furthermore, all elements in $\mathcal{PC}_{n}^{2}(\gamma)$ are obtained by the above construction unrepeatedly.
\end{theorem}

\begin{theorem}\label{thm 6}
	Let $c_{n}(\gamma)$ be a $q$-cyclotomic coset modulo $n$ with $| c_{n}(\gamma) | =\tau$. Assume that $q^{\tau} \equiv 3 \pmod{4}$. Then 
	\begin{description}
		\item[(1)] The cyclotomic $2$-adic integer $\varphi(\gamma)$ gives rise to an element
		$$\mathfrak{c}_{\infty}(\gamma) = (c_{2^{i}n}(\gamma + n\cdot \varphi(\gamma)_{\leq i-1}))_{i \in \mathbb{N}} \in \mathcal{PC}_{n}^{2}(\gamma).$$
		Here $\varphi(\gamma)_{\leq -1}$ is set to be $0$. The sequence $\mathfrak{c}_{\infty}(\gamma)$ is the unique splitting element in $\mathcal{PC}_{n}^{2}(\gamma)$.
		\item[(2)] For each $m \in \mathbb{N}$, the integer $U(\gamma;m)$ gives rise to $2^{v-1}$ elements
		$$\mathfrak{c}_{m,t}(\gamma) = (c_{2^{i}n}(\gamma + n(U(\gamma;m)+2^{m+2}t)))_{i \in \mathbb{N}} \in \mathcal{PC}_{n}^{2}(\gamma),$$
		where $v = v_{2}(q^{\tau}+1)$ and $t \in [0,2^{v-1}-1]$.
	\end{description}
	Furthermore, all elements in $\mathcal{PC}_{n}^{2}(\gamma)$ are obtained by the above construction unrepeatedly.
\end{theorem}

\begin{corollary}\label{coro 2}
	Let the notations be defined as in Theorem \ref{thm 5} and \ref{thm 6}.
	\begin{itemize}
		\item[(1)] If $q^{\tau} \equiv 1 \pmod{4}$, then
		\begin{itemize}
			\item[(1a)] The principal sequence $\mathfrak{c}_{\infty}(\gamma)$ has the same size $\tau$ at every degree.
			\item[(1b)] For any $m \in \mathbb{N}$ and any $t \in [0,2^{v-1}-1]$, the size of each component of $\mathfrak{c}_{m,t}(\gamma)$ is given by
			\begin{small}
				\begin{equation*}
					|c_{2^{i}n}(\gamma + n(U(\gamma;m)+2^{m+1}t)) | = \left\{
					\begin{array}{lcl}
						\tau, \ 0 \leq i \leq m+v;\\
						2^{i-m-v}\tau, \quad i \geq m+v+1.
					\end{array} \right.
				\end{equation*}
			\end{small}
		\end{itemize}
		\item[(2)] If $q^{\tau} \equiv 3 \pmod{4}$, then
		\begin{itemize}
			\item[(2a)] The principal sequence $\mathfrak{c}_{\infty}(\gamma)$ has the same size $\tau = |c_{n/q}(\gamma)|$ at every degree.
			\item[(2b)] For any $m \in \mathbb{N}$ and any $t \in [0,2^{v-1}-1]$, the size of each component of $\mathfrak{c}_{m,\mathbf{t}}$ is given by
			\begin{small}
				\begin{equation*}
					| c_{2^{i}n}(\gamma + n(U(\gamma;m)+2^{m+2}t)) | = \left\{
					\begin{array}{lcl}
						\tau, \ 0 \leq i \leq m+1;\\
						2\tau, \ m+2 \leq i \leq m+v+1;\\
						2^{i-m-v}\tau, \ i \geq m+v+2.
					\end{array} \right.
				\end{equation*}
			\end{small}
		\end{itemize}
	\end{itemize}
\end{corollary}

As in the case that $\ell$ is an odd prime, we call $\mathfrak{c}_{\infty}(\gamma)$ the principal sequence over $c_{n/q}(\gamma)$, and the other elements in $\mathcal{PC}_{n/q}^{2}(\gamma)$ the stable sequences over $c_{n/q}(\gamma)$. For any $m \in \mathbb{N}$, we define
$$\mathcal{PC}_{n}^{2}(\gamma;m) = \{\mathfrak{c}_{m,t}(\gamma) \ | \ t \in [0,2^{v-1}-1]\},$$
where
\begin{equation*}
	v =  \left\{
	\begin{array}{lcl}
		v_{2}(q^{\tau}-1), \ q^{\tau} \equiv 1 \pmod{4};\\
		v_{2}(q^{\tau}+1), \ q^{\tau} \equiv 3 \pmod{4}.
	\end{array} \right.
\end{equation*}
And, further, we define
$$\mathcal{PC}_{n}^{2}(\gamma;\infty) = \{\mathfrak{c}_{\infty}(\gamma)\}.$$
The elements in $\mathcal{PC}_{n}^{2}(\gamma;m)$ are said to be of quasi-stable degree $m$, and of stable degree $m+v$ and $m+v+1$, in the cases where $q^{\tau} \equiv 1 \pmod{4}$ and where $q^{\tau} \equiv 3 \pmod{4}$ respectively.

\begin{corollary}
	The $2$-adic $q$-cyclotomic system $\mathcal{PC}_{n}^{2}(\gamma)$ over $c_{n}(\gamma)$ can be naturally classified via the quasi-stable degree as
	$$\mathcal{PC}_{n}^{2}(\gamma) = \bigsqcup_{m \in \mathbb{N}\cup\{\infty\}}\mathcal{PC}_{n}^{2}(\gamma;m).$$
\end{corollary}

Parallel to Theorem \ref{thm 7}, the structure of the $2$-adic $q$-cyclotomic system $\mathcal{PC}_{n/q}^{2}$ is characterized by the theorem below.

\begin{theorem}
	Let $q$ be a power of an odd prime $p$, and $n$ be a positive odd integer coprime to $q$. Then the $2$-adic $q$-cyclotomic system $\mathcal{PC}_{n}^{2}$ can be decomposed as
	$$\mathcal{PC}_{n}^{2} = \bigsqcup_{\gamma \in \mathcal{CR}_{n}}\left(\bigsqcup_{m \in \mathbb{N}\cup\{\infty\}}\mathcal{PC}_{n}^{2}(\gamma;m)\right).$$
	Furthermore, we have
	\begin{description}
		\item[(1)] The sequences lying in $\mathcal{PC}_{n}^{\ell}(\gamma;\infty)$, $\gamma \in \mathcal{CR}_{n}$, are exactly the sequences that are splitting.
		\item[(2)] Let $\gamma \in \mathcal{CR}_{n}$ with $\tau = c_{n}(\gamma)$. If $q^{\tau} \equiv 1 \pmod{4}$, then for any $m \in \mathbb{N}$, each sequence lying in $\mathcal{PC}_{n}^{2}(\gamma;m)$ is splitting at degree $i \leq m+v-1$, while is stable at degree $i \geq m+v$. If $q^{\tau} \equiv 3 \pmod{4}$, for any $m \in \mathbb{N}$, each sequence in $\mathcal{PC}_{n}^{2}(\gamma;m)$ is splitting at degree $i \leq m$ and at degree $m+2 \leq i \leq m+v$, while is non-splitting at degree $i=m+1$ and at degree $i \geq m+v+1$.
	\end{description}
\end{theorem}

\begin{examples}
	Let $q$ be a power of an odd prime $p$. There is only one $q$-cyclotomic coset modulo $1$, that is, $c_{1}(0)$, therefore one has $\mathcal{PC}_{1}^{2} = \mathcal{PC}_{1}^{2}(0)$. The cyclotomic $2$-adic integer associated to $0$ is 
	$$\varphi(0)=0,$$
	and for any $m \in \mathbb{N}$
	$$U(0;m) = 2^{m}.$$
	
	If $q \equiv 1 \pmod{4}$, set $v = v_{2}(q-1)$. Then $\varphi(0)$ induces 
	$$\mathfrak{c}_{\infty}(0) = (c_{2^{i}}(0))_{i \in \mathbb{N}} \in \mathcal{PC}_{1}^{2},$$
	and for $m \in \mathbb{N}$, $U(0;m)$ induces $2^{v-1}$ sequences
	$$\mathfrak{c}_{m,t}(0) = (c_{2^{i}}(2^{m}(1+2t)))_{i \in \mathbb{N}} \in \mathcal{PC}_{1}^{2},$$
	where $t \in [0,2^{v-1}-1]$. These are all the elements in $\mathcal{PC}_{1}^{2}$. Furthermore, every component $c_{2^{i}}(0)$ of $\mathfrak{c}_{\infty}(0)$ has size $1$, and the size of each component of $\mathfrak{c}_{m,t}(0)$ is given by
	\begin{equation*}
		|c_{2^{i}}(2^{m}(1+2t))| =
		\left\{
		\begin{array}{lcl}
			1, \ \mathrm{if} \ 0 \leq i \leq m+v;\\
			2^{i-m-v}, \ \mathrm{if} \ i \geq m+v+1.
		\end{array} \right.
	\end{equation*}
	
	If $q \equiv 3 \pmod{4}$, set $v = v_{2}(q+1)$. Then $\varphi(0)$ induces 
	$$\mathfrak{c}_{\infty}(0) = (c_{2^{i}}(0))_{i \in \mathbb{N}} \in \mathcal{PC}_{1}^{2},$$
	and for $m \in \mathbb{N}$, $U(0;m)$ induces $2^{v-1}$ sequences
	$$\mathfrak{c}_{m,t}(0) = (c_{2^{i}}(2^{m}(1+4t)))_{i \in \mathbb{N}} \in \mathcal{PC}_{1}^{2},$$
	where $t \in [0,2^{v-1}-1]$. These are all the elements in $\mathcal{PC}_{1}^{2}$. Furthermore, every component $c_{2^{i}}(0)$ of $\mathfrak{c}_{\infty}(0)$ has size $1$, and the size of each component of $\mathfrak{c}_{m,t}(0)$ is given by
	\begin{equation*}
		|c_{2^{i}}(2^{m}(1+4t))| =
		\left\{
		\begin{array}{lcl}
			1, \ \mathrm{if} \ 0 \leq i \leq m+1;\\
			2, \ \mathrm{if} \ m+2 \leq i \leq m+v+1;\\
			2^{i-m-v}, \ \mathrm{if} \ i \geq m+v+2.
		\end{array} \right.
	\end{equation*}
\end{examples}	

\section{Explicit representatives and sizes of cyclotomic cosets}\label{sec 1}
Let $n$ be a positive integer coprime to $q$, and $\lambda$ be a nonzero element in $\mathbb{F}_{q}$. Recall that a $\lambda$-constacyclic code of length $n$ over $\mathbb{F}_{q}$ can be realized as an ideal in the ring $\mathcal{R}_{n,\lambda} = \mathbb{F}_{q}[X]/(X^{n}-\lambda)$, and thus is generated by a unique factor $f(X)$ of $X^{n}-\lambda$. Therefore to determine all $\lambda$-constacyclic codes of length $n$ over $\mathbb{F}_{q}$ needs the irreducible factorization of $X^{n}-\lambda$. Denoting the order of $\lambda$ by $g$, then there is a primitive $ng$-th root $\zeta_{ng}^{h}$ of unity, where $\mathrm{gcd}(h,ng)=1$, such that $\zeta_{ng}^{hn}=\lambda$. Then $X^{n}-\lambda$ can be decomposed as
$$X^{n}-\lambda = (X-\zeta_{ng}^{h})(X-\zeta_{ng}^{h+g})\cdots(X-\zeta_{ng}^{h+(n-1)g})$$
over the splitting field of $X^{n}-\lambda$, which indicates that the defining set of $X^{n}-\lambda$ is
$$\{h,h+g,\cdots,h+(n-1)g\} = h+ g\mathbb{Z}/ng\mathbb{Z}.$$
To determine the irreducible factorization of $X^{n}-\lambda$ is equivalent to giving the partition of $h+ g\mathbb{Z}/ng\mathbb{Z}$ into $q$-cyclotomic cosets modulo $ng$.

On the other hand, since $h$ is coprime to $ng$, the bijection
$$\mu_{h}:1+ g\mathbb{Z}/ng\mathbb{Z} \rightarrow h+ g\mathbb{Z}/ng\mathbb{Z}; \ 1+kg \mapsto h(1+kg)$$
preserves cyclotomic cosets, that is, for any $0 \leq k \leq n-1$
$$\mu_{h}(c_{ng}(1+kg)) = c_{ng}(h(1+kg)).$$
Furthermore, if 
$$1+ g\mathbb{Z}/ng\mathbb{Z} = \bigsqcup_{k}c_{ng}(1+kg)$$
is the partition of $1+ g\mathbb{Z}/ng\mathbb{Z}$ into cyclotomic cosets, then the partition of $h+ g\mathbb{Z}/ng\mathbb{Z}$ is
$$h+ g\mathbb{Z}/ng\mathbb{Z} = \bigsqcup_{k}c_{ng}(h(1+kg)).$$

Motivated by the above facts, we compute the representatives, sizes and enumeration of the $q$-cyclotomic cosets modulo $ng$ contained in $1+ g\mathbb{Z}/ng\mathbb{Z}$. Notice that if $g=1$ then $1+ g\mathbb{Z}/ng\mathbb{Z} = \mathbb{Z}/n\mathbb{Z}$, so determining the full set of representatives, sizes and enumeration of all $q$-cyclotomic cosets modulo $n$ form a special case of this problem.

\subsection{The preimage of a cyclotomic coset along $\widehat{\pi}$}
In this subsection we consider the preimage $\widehat{\pi}_{nn^{\prime}/n}^{-1}(c_{n}(\gamma))$ of a given coset $c_{n}(\gamma)$ along the surjection $\widehat{\pi}_{nn^{\prime}/n}: \mathcal{C}_{nn^{\prime}} \rightarrow \mathcal{C}_{n}$, where $n$, $n^{\prime}$ and $q$ are pairwise coprime. It serves as a preparation for the computation of the partition of $1+g\mathbb{Z}/ng\mathbb{Z}$ into $q$-cyclotomic cosets, which is given in Subsection. 

First we deal with the subcase where $n^{\prime}=\ell^{j}$ is a power of a prime number $\ell$ which is different $p$ and does not divide $n$. For any $\mathfrak{c} = (c_{\ell^{i}n}(\gamma_{i}))_{i \in \mathbb{N}} \in \mathcal{PC}_{n}^{\ell}(\gamma)$ and any positive integer $j$, we denote by
$$\mathfrak{c}_{j} = c_{\ell^{j}n}(\gamma_{j})$$
the component of degree $j$ of $\mathfrak{c}$. And for each $m \in \mathbb{N}$, we write
$$\mathcal{PC}_{n}^{\ell}(\gamma;m)_{j} = \{\mathfrak{c}_{j} \ | \ \mathfrak{c} \in \mathcal{PC}_{n}^{\ell}(\gamma;m)\}.$$

\begin{lemma}\label{lem 5}
	For any $j \in \mathbb{N}^{+}$, the preimage $\widehat{\pi}_{\ell^{j}n/n}^{-1}(c_{n}(\gamma))$ along $\widehat{\pi}_{\ell^{j}n/n}$ is
	\begin{equation}\label{eq 12}
		\widehat{\pi}_{\ell^{j}n/n}^{-1}(c_{n}(\gamma)) = \bigsqcup_{m \in [0,j-1]\cup\{\infty\}}\mathcal{PC}_{n}^{\ell}(\gamma;m)_{j}.
	\end{equation}
\end{lemma}

\begin{proof}
	The cyclotomic cosets contained in the preimage $\widehat{\pi}_{\ell^{j}n/n}^{-1}(c_{n}(\gamma))$ can be obtained by taking the degree-$j$ components of the sequences in $\mathcal{PC}_{n}^{\ell}(\gamma;m)$. Now Lemma \ref{lem 5} follows from the fact that the sequences lying in $\mathcal{PC}_{n}^{\ell}(\gamma;m)$ for $m \geq j$ and $m=\infty$ have the same degree-$j$ component.
\end{proof}

Combining the decomposition \eqref{eq 12} and the explicit expression of $\mathcal{PC}_{n}^{\ell}(\gamma;m)$, we can obtain the cyclotomic cosets which consist $\widehat{\pi}_{\ell^{j}n/n}^{-1}(c_{n}(\gamma))$ concretely. In order to save notations we apply the following conventions. Denote by $\tau$ the size of $c_{n}(\gamma)$, and by $r$ the index of the subgroup $\langle q^{\tau}\rangle$ in $(\mathbb{Z}/n\mathbb{Z})^{\ast}$. Set
\begin{equation*}
	v =\left\{
	\begin{array}{lcl}
		v_{2}(q^{\tau}+1), \ \mathrm{if} \ \ell = 2 \ \mathrm{and} \ q^{\tau} \equiv 3 \pmod{4};\\
		v_{\ell}(q^{\frac{\tau(\ell-1)}{r}}-1), \quad \mathrm{otherwise}.
	\end{array} \right.
\end{equation*}
If $\ell=2$, $U_{n}^{\ell}(\gamma;m,1)$ just stands for $U_{n}^{2}(\gamma;m)$.

\begin{proposition}\label{prop 3}
	Let the notations be defined as above. 
	\begin{description}
		\item[(1)] If either $\ell$ is odd, or $\ell=2$ and $q^{\tau} \equiv 1 \pmod{4}$, then we have
		\begin{small}
		\begin{align*}
			\widehat{\pi}_{\ell^{j}n/n}^{-1}(c_{n}(\gamma)) =& \{c_{\ell^{j}n}(\gamma+n\cdot \varphi_{n}^{\ell}(\gamma)_{\leq j-1})\} \sqcup\\ &\left(\bigsqcup_{m=0}^{j-1}\{c_{\ell^{j}n}(\gamma+n(U_{n}^{\ell}(\gamma;m,a)+\ell^{m+1}t)) \ | \ a \in [1,r], t \in [0,\ell^{\mathrm{min}\{v,j-m\}-1}-1]\}\right).
		\end{align*}	
		\end{small}
		\item[(2)] If $\ell=2$ and $q^{\tau} \equiv 3 \pmod{4}$, then
		\begin{small}
			\begin{align*}
				\widehat{\pi}_{2^{j}n/n}^{-1}(c_{n}(\gamma)) =& \{c_{2^{j}n}(\gamma+n\cdot \varphi_{n}^{2}(\gamma)_{\leq j-1})\} \sqcup\\ &\left(\bigsqcup_{m=0}^{j-1}\{c_{2^{j}n}(\gamma+n(U_{n}^{2}(\gamma;m)+2^{m+2}t)) \ | \ t \in [0,2^{\mathrm{max}\{0,\mathrm{min}\{v-1,j-m-2\}\}}-1]\}\right).
			\end{align*}	
		\end{small}
	\end{description}
\end{proposition}

\begin{proof}
	We will prove for the case where $\ell$ is odd and $\ell \mid q^{\tau}-1$. The proofs for the other cases are similar. By Lemma \ref{lem 5} it suffices to compute $\mathcal{PC}_{n}^{\ell}(\gamma;m)_{j}$ for $m=\infty$ and $0 \leq m \leq j-1$.

	For $m=\infty$, $\mathcal{PC}_{n}^{\ell}(\gamma;m) = \{\mathfrak{c}_{\infty}(\gamma)\}$, and thus
	\begin{equation}
		\mathcal{PC}_{n}^{\ell}(\gamma;\infty)_{j} = \{c_{\ell^{j}n}(\gamma+n\cdot\varphi_{n}^{\ell}(\gamma)_{\leq j-1})\}.
	\end{equation}
	
	Let $m \in [0,j-1]$. Suppose first that $m$ is such that $m+v \leq j$. Then for any $\mathfrak{c}_{m,a,t}(\gamma) \in \mathcal{PC}_{n}^{\ell}(\gamma;m)$, its degree-$j$ component is $c_{\ell^{j}n}(\gamma+n(U_{n}^{\ell}(\gamma;m,a)+\ell^{m+1}t))$. By Theorem \ref{thm 7} every $\mathfrak{c}_{m,a,t}(\gamma)$ is stable at degree $i \geq j$, therefore the cosets $c_{\ell^{j}n}(\gamma+n(U_{n}^{\ell}(\gamma;m,a)+\ell^{m+1}t))$, $a \in [1,\ell-1]$ and $t \in [0,\ell^{v-1}-1]$, are pairwise distinct. Hence we have
	\begin{equation}
		\mathcal{PC}_{n}^{\ell}(\gamma;m)_{j} = \{c_{\ell^{j}n}(\gamma+n(U_{n}^{\ell}(\gamma;m,a)+\ell^{m+1}t)) \ | \ a \in [1,\ell-1], t \in [0,\ell^{\mathrm{min}\{v,j-m\}-1}-1]\}.
	\end{equation}
	
	Next suppose that $m+v > j$. Since $m \leq j-1$, any $\mathfrak{c}_{m,a,t}(\gamma)$ lying in $\mathcal{PC}_{n}^{\ell}(\gamma;m)$ has component at degree $m+1$ being
	$$c_{\ell^{m}n}(\gamma+n(U_{n}^{\ell}(\gamma;m,a)+\ell^{m+1}t)) = c_{\ell^{m}n}(\gamma+n\cdot U_{n}^{\ell}(\gamma;m,a)).$$
	As $v_{\ell}(\gamma+n\cdot U_{n}^{\ell}(\gamma;m,a))=m$, applying Lemma \ref{lem 1} successively yields
	\begin{align*}
		\widehat{\pi}_{\ell^{j}n/\ell^{m}n}^{-1}(c_{n}(\gamma)) &= \{c_{\ell^{j}n}(\gamma+n(U_{n}^{\ell}(\gamma;m,a)+\sum_{k=1}^{j-m-1}t_{k}\ell^{m+k}))\ | \ (t_{1},\cdots,t_{j-m-1}) \in [0,\ell-1]^{j-m-1}\}\\
		&= \{c_{\ell^{j}n}(\gamma+n(U_{n}^{\ell}(\gamma;m,a)+\ell^{m+1}t)) \ | \ t \in [0,\ell^{j-m-1}-1]\}.
	\end{align*}
	It follows that
	\begin{equation}\label{eq 11}
		\mathcal{PC}_{n}^{\ell}(\gamma;m)_{j} = \{c_{\ell^{j}n}(\gamma+n(U_{n}^{\ell}(\gamma;m,a)+\ell^{m+1}t)) \ | \ a \in [1,\ell-1], t \in [0,\ell^{j-m-1}-1]\}.
	\end{equation}
	
	Now the equalities \eqref{eq 12}--\eqref{eq 11} together imply the conclusion.
\end{proof}

\begin{corollary}\label{coro 5}
	\begin{description}
		\item[(1)] If either $\ell$ is odd, or $\ell=2$ and $q^{\tau} \equiv 1 \pmod{4}$, then the sizes of the cosets in $\widehat{\pi}_{\ell^{j}n/n}^{-1}(c_{n}(\gamma))$ are given by
		$$|c_{\ell^{j}n}(\gamma+n\cdot \varphi_{n}^{\ell}(\gamma)_{\leq j-1})| = \tau,$$
		and
		$$|c_{\ell^{j}n}(\gamma+n(U_{n}^{\ell}(\gamma;m,a)+\ell^{m+1}t))| = \ell^{\mathrm{max}\{0,j-m-v\}}\cdot \dfrac{\ell-1}{r}\cdot \tau$$
		for all $a \in [1,r]$ and $t \in [0,\ell^{\mathrm{min}\{v,j-m\}-1}-1]$.
		\item[(2)] If $\ell=2$ and $q^{\tau} \equiv 3 \pmod{4}$, then the sizes of the cosets in $\widehat{\pi}_{\ell^{j}n/n}^{-1}(c_{n}(\gamma))$ are given by
		$$|c_{2^{j}n}(\gamma+n\cdot \varphi_{n}^{2}(\gamma)_{\leq j-1})| = \tau,$$
		and
		\begin{equation*}
			|c_{2^{j}n}(\gamma+n(U_{n}^{2}(\gamma;m)+2^{m+2}t))| =
			\left\{
			\begin{array}{lcl}
				\tau, \ \mathrm{if} \ m=j-1;\\
				2^{\mathrm{max}\{1,j-m-v\}}\cdot\tau, \ \mathrm{if} \ 0 \leq m \leq j-2
			\end{array} \right.
		\end{equation*}
		for all $t \in [0,2^{\mathrm{min}\{v,j-m-1\}-1}-1]$.
	\end{description}
\end{corollary}

\begin{proof}
	The conclusions follow from Corollary \ref{coro 1}, \ref{coro 3}, \ref{coro 2} and Proposition \ref{prop 3} directly.
\end{proof}

Further, the enumeration of $q$-cyclotomic cosets modulo $\ell^{j}n$ contained in $\widehat{\pi}_{\ell^{j}n/n}^{-1}(c_{n}(\gamma))$ can be obtained as follow.

\begin{corollary}\label{coro 6}
	\begin{description}
		\item[(1)] If either $\ell$ is odd, or $\ell=2$ and $q^{\tau} \equiv 1 \pmod{4}$, then there are in total 
		$$1+\dfrac{r(\ell^{\mathrm{min}\{j,v\}}-1)}{\ell-1}+\mathrm{max}\{0,j-v\}\cdot r\ell^{v-1}$$
		$q$-cyclotomic cosets modulo $\ell^{j}n$ lying in $\widehat{\pi}_{\ell^{j}n/n}^{-1}(c_{n}(\gamma))$.
		\item[(2)] If $\ell=2$ and $q^{\tau} \equiv 3 \pmod{4}$, then there are in total
		$$1+2^{\mathrm{min}\{j-1,v\}}+\mathrm{max}\{0,j-v-1\}\cdot 2^{v-1}$$
		$q$-cyclotomic cosets modulo $2^{j}n$ lying in $\widehat{\pi}_{2^{j}n/n}^{-1}(c_{n}(\gamma))$.
	\end{description}
\end{corollary}

\begin{proof}
	Here we assume that either $\ell$ is odd, or $\ell=2$ and $q^{\tau} \equiv 1 \pmod{4}$. The result for the latter case can be proved similarly. By Proposition \ref{prop 3} there are
	$$1+\sum_{m=0}^{j-1}r\ell^{\mathrm{min}\{v,j-m\}-1} = 1+ r\sum_{m=0}^{j-1}\ell^{\mathrm{min}\{v,j-m\}-1}$$
	cosets in $\widehat{\pi}_{\ell^{j}n/n}^{-1}(c_{n}(\gamma))$. If $j \leq v$ then
	$$\sum_{m=0}^{j-1}\ell^{\mathrm{min}\{v,j-m\}-1} = \sum_{k=0}^{j-1}\ell^{k} = \dfrac{\ell^{j}-1}{\ell-1};$$
	while if $j >v$ then
	$$\sum_{m=0}^{j-1}\ell^{\mathrm{min}\{v,j-m\}-1} = \sum_{k=0}^{v-1}\ell^{k} + (j-v)\cdot\ell^{v-1} = \dfrac{\ell^{v}-1}{\ell-1} + (j-v)\ell^{v-1}.$$
	Now the conclusion follows immediately.
\end{proof}

Now we turn to the general case. Assume that $n^{\prime}=p_{1}^{j_{1}}\cdots p_{s}^{j_{s}}$, where $p_{1} < \cdots < p_{s}$ are prime numbers that are different from $p$ and do not divide $n$, and $j_{1},\cdots,j_{s}$ are positive integers. Denote by $\overline{[0,j_{k}-1]} = [0,j_{k}-1] \cup \{\infty\}$ for any $1 \leq k \leq s$, and let $\mathfrak{m} = (m_{1},\cdots,m_{s})$ be a $s$-tuple of integers $m_{k} \in \overline{[0,j_{k}-1]}$. We introduce the following notations.

\begin{itemize}
	\item[1.] For any $1\leq k\leq s$, $\mathfrak{m}_{\leq k}$ denotes the $k$-tuple $(m_{1},\cdots,m_{k})$.
	\item[2.] For $k=1$, we set $r_{1} = [(\mathbb{Z}/p_{1}\mathbb{Z})^{\ast}: \langle q^{\tau}\rangle]$ and
	\begin{equation*}
		v_{1} =
		\left\{
		\begin{array}{lcl}
			v_{2}(q^{\tau}+1), \ \mathrm{if} \ p_{1}=2 \ \mathrm{and} \ q^{\tau} \equiv 3 \pmod{4};\\
			v_{p_{1}}(q^{\frac{\tau(p_{1}-1)}{r_{1}}}-1), \ \mathrm{otherwise}.
		\end{array} \right.
	\end{equation*}
	If $p_{1}=2$ and $q^{\tau} \equiv 3 \pmod{4}$, set
	\begin{equation*}
		\tau_{1}(\mathfrak{m}_{\leq 1}) =
		\left\{
		\begin{array}{lcl}
			\tau, \ \mathrm{if} \ m_{1}=j_{1}-1 \ \mathrm{or} \ \infty;\\
			2^{\mathrm{max}\{1,j_{1}-m_{1}-v_{1}\}}\cdot\tau, \ \mathrm{if} \ 0 \leq m_{1} \leq j_{1}-2.
		\end{array} \right.
	\end{equation*}
	Otherwise, set
	\begin{equation*}
		\tau_{1}(\mathfrak{m}_{\leq 1}) =
		\left\{
		\begin{array}{lcl}
			\tau, \ \mathrm{if} \ m_{1}=\infty;\\
			p_{1}^{\mathrm{max}\{0,j_{1}-m_{1}-v_{1}\}}\cdot\dfrac{p_{1}-1}{r_{1}}\cdot\tau, \ \mathrm{if} \ 0 \leq m_{1} \leq j_{1}-1.
		\end{array} \right.
	\end{equation*}
	\item[3.] For $2 \leq k \leq s$, we set $r_{k}(\mathfrak{m}_{\leq k-1})= [(\mathbb{Z}/p_{k}\mathbb{Z})^{\ast}: \langle q^{\widetilde{\tau}_{k}}\rangle]$, $v_{k}(\mathfrak{m}_{\leq k-1}) = v_{p_{k}}(q^{\widetilde{\tau}_{k}}-1)$, and
	\begin{equation*}
		\tau_{k}(\mathfrak{m}_{\leq k}) =
		\left\{
		\begin{array}{lcl}
			\tau_{k-1}(\mathfrak{m}_{\leq k-1}), \ \mathrm{if} \ m_{k}=\infty;\\
			p_{k}^{\mathrm{max}\{0,j_{k}-m_{k}-v_{k}(\mathfrak{m}_{\leq k-1})\}}\cdot\dfrac{p_{k}-1}{r_{k}(\mathfrak{m}_{\leq k-1})}\cdot\tau, \ \mathrm{if} \ 0 \leq m_{k} \leq j_{k}-1,
		\end{array} \right.
	\end{equation*}
	where $\widetilde{\tau}_{k} = \frac{\tau_{k-1}(\mathfrak{m}_{\leq k-1})(p_{k}-1)}{r_{k}(\mathfrak{m}_{\leq k-1})}$.
\end{itemize}

\begin{definition}
	Given any $\mathfrak{m} = (m_{1},\cdots,m_{s})$, where $m_{k} \in \overline{[0,j_{k}-1]}$ for each $1 \leq k \leq s$, define a set $\mathcal{I}_{nn^{\prime}/n}(\gamma;\mathfrak{m})$ by
	\begin{align*}
		\mathcal{I}_{nn^{\prime}/n}(\gamma;\mathfrak{m}) = \{(\mathfrak{a},\mathfrak{t}) \ | \ & \mathfrak{a} = (a_{1},\cdots,a_{s}) \in \mathbb{N}^{s} \ \mathrm{and} \ \mathfrak{t} = (t_{1},\cdots,t_{s}) \in \mathbb{N}^{s} \\ 
		&\mathrm{satisfy} \ \mathrm{Condition} \ (\romannumeral1)--(\romannumeral3)\},
	\end{align*}
	where the conditions (\romannumeral1)--(\romannumeral3) are listed as follows:
	\begin{description}
		\item[(\romannumeral1)] For every $1 \leq k \leq s$, if $m_{k}=\infty$ then $a_{k}=t_{k}=0$.
		\item[(\romannumeral2)] Assume that $m_{1} \in [0,j_{1}-1]$. If $p_{1}=2$ and $q^{\tau}\equiv 3 \pmod{4}$ then $a_{1}=1$ and $t_{1} \in [0,2^{\mathrm{max}\{0,\mathrm{min}\{v-1,j-m-2\}\}}-1]$; otherwise $a_{1} \in [1,r_{1}]$ and $t_{1} \in [0,p_{1}^{\mathrm{min}\{v_{1},j_{1}-m_{1}\}-1}-1]$.
		\item[(\romannumeral3)] For every $2 \leq k \leq s$, if $m_{k} \in [0,j_{k}-1]$ then $a_{k} \in [1,r_{k}(\mathfrak{m}_{\leq k-1})]$ and $t_{k} \in [0,p_{k}^{\mathrm{min}\{v_{k}(\mathfrak{m}_{\leq k-1}),j_{k}-m_{k}\}-1}-1]$.
	\end{description}
	Moreover, we define
	$$\mathcal{I}_{nn^{\prime}/n}(\gamma) = \bigsqcup_{\mathfrak{m} \in \overline{[0,j_{1}-1]}\times\cdots\times\overline{[0,j_{s}-1]}}\mathcal{I}_{nn^{\prime}/n}(\gamma;\mathfrak{m}).$$
\end{definition}

\begin{theorem}\label{thm 1}
	There is an one-to-one correspondence from $\mathcal{I}_{nn^{\prime}/n}(\gamma)$ onto $\widehat{\pi}_{nn^{\prime}/n}^{-1}(c_{n}(\gamma))$.
\end{theorem}

\begin{proof}
	We construct a map
	$$\omega: \mathcal{I}_{nn^{\prime}/n}(\gamma) \rightarrow \widehat{\pi}_{nn^{\prime}/n}^{-1}(c_{n}(\gamma)); \ I \mapsto c_{nn^{\prime}}(\gamma_{I})$$
	as follow. Each element $I$ in $\mathcal{I}_{nn^{\prime}/n}(\gamma)$ can be presented as a $3$-tuple $I = (\mathfrak{m},\mathfrak{a},\mathfrak{t})$, where $\mathfrak{m} \in \overline{[0,j_{1}-1]}\times\cdots\times\overline{[0,j_{s}-1]}$ and $(\mathfrak{a},\mathfrak{t}) \in \mathcal{I}_{nn^{\prime}/n}(\gamma;\mathfrak{m})$. We define the representative $\gamma_{I}$ inductively by
	\begin{itemize}
		\item[1.] For $k=1$,
		\begin{equation*}
			\gamma_{I,1} =
			\left\{
			\begin{array}{lcl}
				\gamma+n\cdot\varphi_{n}^{p_{1}}(\gamma)_{\leq j_{1}-1}, \ \mathrm{if} \ m_{1}=\infty;\\
				\gamma+n(U_{n}^{2}(\gamma;m_{1})+2^{m_{1}+2}t_{1}), \ \mathrm{if} \ 0 \leq m_{1} \leq j_{1}-1, \ p_{1}=2, \ \mathrm{and} \ q^{\tau} \equiv 1 \pmod{4};\\
				\gamma+n(U_{n}^{p_{1}}(\gamma;m_{1},a_{1})+p_{1}^{m_{1}+1}t_{1}), \ \mathrm{otherwise}.
			\end{array} \right.
		\end{equation*}
		\item[2.] For $2 \leq k \leq s$,
		\begin{equation*}
			\gamma_{I,k} =
			\left\{
			\begin{array}{lcl}
				\gamma_{I,k-1}+p_{1}^{j_{1}}\cdots p_{k-1}^{j_{k-1}}n\cdot\varphi_{p_{1}^{j_{1}}\cdots p_{k-1}^{j_{k-1}}n}^{p_{k}}(\gamma_{I,k-1})_{\leq j_{k}-1}, \ \mathrm{if} \ m_{k}=\infty;\\
				\gamma_{I,k-1}+p_{1}^{j_{1}}\cdots p_{k-1}^{j_{k-1}}n(U_{p_{1}^{j_{1}}\cdots p_{k-1}^{j_{k-1}}n}^{p_{k}}(\gamma_{I,k-1};m_{k},a_{k})+p_{1}^{m_{k}+1}t_{k}), \ \mathrm{if} \ 0 \leq m_{k} \leq j_{k}-1.
			\end{array} \right.
		\end{equation*}
		\item[3.] $\gamma_{I} = \gamma_{I,s}$.
	\end{itemize}
	Applying Proposition \ref{prop 3} and Corollary \ref{coro 5} successively yields that $\omega$ is a bijection.
\end{proof}

\begin{corollary}
	Let the notations be defined as in Theorem \ref{thm 1}. Then for any $I=(\mathfrak{m},\mathfrak{a},\mathfrak{t}) \in \mathcal{I}_{nn^{\prime}/n}(\gamma)$ the size of $c_{nn^{\prime}}(\gamma_{I})$ is 
	$$|c_{nn^{\prime}}(\gamma_{I})| = \tau_{s}(\mathfrak{m}).$$
\end{corollary}

\subsection{Explicit representatives and sizes of the cyclotomic cosets in $1+g\mathbb{Z}/ng\mathbb{Z}$}\label{sec 5}
Let $g$ be a positive divisor of $q-1$, and $n$ be a positive integer coprime to $q$. As the motivation being explained at the beginning of Section \ref{sec 1}, we determine the $q$-cyclotomic cosets modulo $ng$ that are contained in the set $1+g\mathbb{Z}/ng\mathbb{Z}$ in this subsection, which leads to the partition of $1+g\mathbb{Z}/ng\mathbb{Z}$ into cyclotomic cosets. Notice that $c_{g}(1) = \{1\}$ as $g$ divides $q-1$. It follows that $1+g\mathbb{Z}/ng\mathbb{Z}$ is exactly the preimage of $c_{g}(1)$ along the canonical projection $\pi_{ng/g}$.

Assume that 
$$g = g_{1}^{e_{1}}\cdots g_{u}^{e_{u}}$$
and
$$n = g_{1}^{f_{1}}\cdots g_{u}^{f_{u}}p_{1}^{j_{1}}\cdots p_{s}^{j_{s}},$$
where $g_{1},\cdots,g_{u},p_{1},\cdots,p_{s}$ are pairwise distinct prime numbers that are different from $p$, $e_{1},\cdots,e_{u},j_{1},\cdots,j_{s}$ are positive integers, and $f_{1},\cdots,f_{u}$ are nonnegative integers. Further, we assume that $g_{1}<\cdots<g_{u}$ and $p_{1}<\cdots<p_{s}$. Denote the divisor $g_{1}^{f_{1}}\cdots g_{u}^{f_{u}}$ of $n$ by $n_{1}$. We set
\begin{itemize}
	\item[1.] $y_{1} = \mathrm{max}\{0,\mathrm{min}\{f_{1}-1,v_{g_{1}(q+1)-e_{1}}\}\}$ if $g_{1}=2$ and $q \equiv 3 \pmod{4}$; otherwise $y_{1} = \mathrm{min}\{f_{1},v_{g_{1}}(q-1)-e_{1}\}$, and
	\item[2.] $y_{k} = \mathrm{min}\{f_{k},v_{g_{k}}(q-1)-e_{k}\}$ for all $2 \leq k \leq u$.
\end{itemize}

Let $\mathbf{x}$ be a $u$-tuple $\mathbf{x} = (x_{1},\cdots,x_{u})$ where $x_{k} \in [0,g_{k}^{y_{k}}-1]$ for each $1 \leq k \leq u$. We attach to $\mathbf{x}$ an integer
\begin{equation*}
	\gamma_{\mathbf{x}} =
	\left\{
	\begin{array}{lcl}
		1+g(g_{1}x_{1}+g_{1}^{f_{1}}x_{2}+\cdots+g_{1}^{f_{1}}\cdots g_{u-1}^{f_{u-1}}x_{u}), \ \mathrm{if} \ g_{1}=2 \ \mathrm{and} \ q \equiv 3 \pmod{4};\\
		1+g(x_{1}+g_{1}^{f_{1}}x_{2}+\cdots+g_{1}^{f_{1}}\cdots g_{u-1}^{f_{u-1}}x_{u}), \ \mathrm{otherwise}.
	\end{array} \right.
\end{equation*}

\begin{proposition}
	With the notations as above, then we have
	$$\widehat{\pi}_{n_{1}g/g}^{-1}(c_{r}(1)) = \{c_{n_{1}g}(\gamma_{\mathbf{x}}) \ | \ \mathbf{x} \in [0,g_{1}^{y_{1}}-1]\times\cdots[0,g_{u}^{y_{u}}-1]\}.$$
	Furthermore, all the cosets $c_{n_{1}g}(\gamma_{\mathbf{x}})$ have the same size, given by
	\begin{equation*}
		\tau =
		\left\{
		\begin{array}{lcl}
			g_{1}^{\mathrm{max}\{e_{1}+f_{1}-v_{g_{1}}(q+1),1\}}\cdot\prod\limits_{k=2}^{u}g_{k}^{\mathrm{max}\{e_{k}+f_{k}-v_{g_{k}}(q-1),0\}}, \ \mathrm{if} \ g_{1}=2 \ \mathrm{and} \ q \equiv 3 \pmod{4};\\
			\prod\limits_{k=1}^{u}g_{k}^{\mathrm{max}\{e_{k}+f_{k}-v_{g_{k}}(q-1),0\}}, \ \mathrm{otherwise}.
		\end{array} \right.
	\end{equation*}
\end{proposition}

\begin{proof}
	The conclusion follows from successive applications of Lemma \ref{lem 4} and \ref{lem 1}.
\end{proof}

Given any $c_{n_{1}g}(\gamma_{\mathbf{x}})$, where $\mathbf{x} \in [0,g_{1}^{y_{1}}-1]\times\cdots[0,g_{u}^{y_{u}}-1]$, by Theorem \ref{thm 1} the preimage $\widehat{\pi}_{ng/n_{1}g}^{-1}(c_{n_{1}g}(\gamma_{\mathbf{x}}))$ is given by
$$\widehat{\pi}_{ng/n_{1}g}^{-1}(c_{n_{1}g}(\gamma_{\mathbf{x}})) = \{c_{ng}(\gamma_{\mathbf{x},I}) \ | \ I \in \mathcal{I}_{ng/n_{1}g}(\gamma_{\mathbf{x}})\},$$
where $\gamma_{\mathbf{x},I} = (\gamma_{\mathbf{x}})_{I}$ is constructed as in the proof of Theorem \ref{thm 1}. As $\pi_{ng/g} = \pi_{ng/n_{1}g} \circ \pi_{n_{1}g/g}$, one has 
\begin{equation}\label{eq 13}
	\begin{aligned}
		\widehat{\pi}_{n_{1}g/g}^{-1}(c_{r}(1)) &= \widehat{\pi}_{ng/n_{1}g}^{-1}(\widehat{\pi}_{n_{1}g/g}^{-1}(c_{g}(1)))\\
		&= \bigsqcup_{\mathbf{x} \in [0,g_{1}^{y_{1}}-1]\times\cdots[0,g_{u}^{y_{u}}-1]}\widehat{\pi}_{ng/n_{1}g}^{-1}(c_{n_{1}g}(\gamma_{\mathbf{x}}))\\
		&= \bigsqcup_{\mathbf{x} \in [0,g_{1}^{y_{1}}-1]\times\cdots[0,g_{u}^{y_{u}}-1]}\bigsqcup_{I \in \mathcal{I}_{ng/n_{1}g}(\gamma_{\mathbf{x}})}\{c_{ng}(\gamma_{\mathbf{x},I})\}.
	\end{aligned}
\end{equation}

\begin{theorem}\label{thm 8}
	The set $1+g\mathbb{Z}/ng\mathbb{Z}$ can be written as
	$$1+g\mathbb{Z}/ng\mathbb{Z} = \bigsqcup_{\mathbf{x} \in [0,g_{1}^{y_{1}}-1]\times\cdots[0,g_{u}^{y_{u}}-1]}\bigsqcup_{I \in \mathcal{I}_{ng/n_{1}g}(\gamma_{\mathbf{x}})}c_{ng}(\gamma_{\mathbf{x},I}).$$
\end{theorem}

\begin{proof}
	As $1+g\mathbb{Z}/ng\mathbb{Z} = \pi_{ng/g}^{-1}(c_{g}(1))$, then \eqref{eq 13} indicates the result.
\end{proof}

\begin{corollary}
	Let the notations be given as in Theorem \ref{thm 1}. Then for any $\mathbf{x} \in [0,g_{1}^{y_{1}}-1]\times\cdots[0,g_{u}^{y_{u}}-1]$ and $I \in \mathcal{I}_{ng/n_{1}g}(\gamma_{\mathbf{x}})$, the size of the coset $c_{ng}(\gamma_{\mathbf{x},I})$ is $\tau_{s}(\mathfrak{m})$, where $I = (\mathfrak{m},\mathfrak{a},\mathfrak{t})$.
\end{corollary}

\section{Classification of constacyclic codes over finite fields}\label{sec 6}
Based on the results obtained in the last section, we determine all the $\lambda$-constacyclic codes of length $N$ over $\mathbb{F}_{q}$, where $\lambda$ is any nonzero element in $\mathbb{F}_{q}$ and $N$ is any positive integer. Here in order to include the case of repeated-root constacyclic codes we do not require $N$ to be coprime to $q$.

Assume that $N = p^{v}n$ where $v = v_{p}(N) \geq 0$ and $\mathrm{gcd}(n,p)=1$. Denoting by $g$ the order of $\lambda$, there is a unique integer $1 \leq b \leq g-1$ satisfying $bp^{v} \equiv 1 \pmod{g}$. Setting $\lambda^{\prime} = \lambda^{b}$, then one has $(\lambda^{\prime})^{p^{v}} = \lambda$, and
$$X^{N}-\lambda = (X^{n}-\lambda^{\prime})^{p^{v}}.$$
Moreover, the order of $\lambda^{\prime}$ is also $g$. There exists a primitive $ng$-th root of unity $\zeta_{ng}^{h}$, where $h$ is coprime to $ng$, such that $\zeta_{ng}^{hn}=\lambda^{\prime}$. As being explained at the beginning of Section \ref{sec 1}, the defining set of $X^{n}-\lambda^{\prime}$ is
$$h+g\mathbb{Z}/ng\mathbb{Z} = \{h,h+g,\cdots,h+(n-1)g\}.$$

Assume that 
$$g = g_{1}^{e_{1}}\cdots g_{u}^{e_{u}}$$
and
$$n = g_{1}^{f_{1}}\cdots g_{u}^{f_{u}}p_{1}^{j_{1}}\cdots p_{s}^{j_{s}},$$
where $g_{1},\cdots,g_{u},p_{1},\cdots,p_{s}$ are pairwise distinct prime numbers that are different from $p$, $e_{1},\cdots,e_{u},j_{1},\cdots,j_{s}$ are positive integers, and $f_{1},\cdots,f_{u}$ are nonnegative integers. And other notations and conventions in Subsection \ref{sec 5} are also adopted here. Combining Theorem \ref{thm 8} and Lemma \ref{lem 7}, we obtain the irreducible factorization of $X^{N}-\lambda$ over $\mathbb{F}_{q}$. Given $\mathfrak{m} = (m_{1},\cdots,m_{s}) \in \overline{[0,j_{1}-1]}\times\cdots\times\overline{[0,j_{s}-1]}$, we set $z_{k} = \mathrm{min}\{j_{k}-m_{k},1\}$ for every $1 \leq k \leq s$,
\begin{equation*}
	\omega_{\mathfrak{m}} =
	\left\{
	\begin{array}{lcl}
		2\mathrm{ord}_{\mathrm{rad}(g)p_{1}^{z_{1}}\cdots p_{s}^{z_{s}}}(q), \ \mathrm{if} \ q^{\mathrm{ord}_{\mathrm{rad}(g)p_{1}^{z_{1}}\cdots p_{s}^{z_{s}}}(q)} \equiv 3 \pmod{4} \ \mathrm{and} \ 8 \mid n_{1}gp_{1}^{j_{1}-m_{1}}\cdots p_{s}^{j_{s}-m_{s}};\\
		\mathrm{ord}_{\mathrm{rad}(g)p_{1}^{z_{1}}\cdots p_{s}^{z_{s}}}(q), \ \mathrm{otherwise},
	\end{array} \right.
\end{equation*}
and $d_{\mathfrak{m}} = \frac{n\omega_{\mathfrak{m}}}{\tau_{s}(\mathfrak{m})}$.

\begin{proposition}
	The irreducible factorization of $X^{N}-\lambda$ over $\mathbb{F}_{q}$ is given by
	\begin{equation}\label{eq 14}
		X^{N}-\lambda = \prod_{\mathbf{x} \in [0,g_{1}^{y_{1}}-1]\times\cdots[0,g_{u}^{y_{u}}-1]}\prod_{I \in \mathcal{I}_{ng/n_{1}g}(\gamma_{\mathbf{x}})}M_{h\gamma_{\mathbf{x},I}}(X)^{p^{v}}.
	\end{equation}
	Moreover, for any $\mathbf{x} \in [0,g_{1}^{y_{1}}-1]\times\cdots[0,g_{u}^{y_{u}}-1]$ and $I = (\mathfrak{m},\mathfrak{a},\mathfrak{t}) \in \mathcal{I}_{ng/n_{1}g}(\gamma_{\mathbf{x}})$, the irreducible polynomial $M_{h\gamma_{\mathbf{x},I}}(X)$ can be expressed as
	$$M_{h\gamma_{\mathbf{x},I}}(X) = \sum_{k=0}^{\omega_{\mathfrak{m}}}(-1)^{\omega_{\mathfrak{m}}-k}\sigma_{\omega_{\mathfrak{m}}-k}(\zeta_{d_{\mathfrak{m}}}^{h\gamma_{\mathbf{x},I}},\zeta_{d_{\mathfrak{m}}}^{h\gamma_{\mathbf{x},I} q},\cdots,\zeta_{d_{\mathfrak{m}}}^{h\gamma_{\mathbf{x},I} q^{\omega_{\mathfrak{m}}-1}})X^{\frac{\tau_{s}(\mathfrak{m}) k}{\omega_{\mathfrak{m}}}}.$$
\end{proposition}

\begin{proof}
	Recall that the map
	$$\mu_{h}:1+ g\mathbb{Z}/ng\mathbb{Z} \rightarrow h+ g\mathbb{Z}/ng\mathbb{Z}; \ 1+kg \mapsto h(1+kg)$$
	is a bijection that preserves cyclotomic cosets, therefore without losing generality we may suppose that $h=1$.
	
	The irreducible factorization \eqref{eq 14} follows from Theorem \ref{thm 8}. For any $\mathbf{x} \in [0,g_{1}^{y_{1}}-1]\times\cdots[0,g_{u}^{y_{u}}-1]$ and $I = (\mathfrak{m},\mathfrak{a},\mathfrak{t}) \in \mathcal{I}_{ng/n_{1}g}(\gamma_{\mathbf{x}})$, by the construction of $\gamma_{\mathbf{x},I}$ we have
	$$\mathrm{gcd}(\gamma_{\mathbf{x},I},ng) = p_{1}^{m_{1}}\cdots p_{s}^{m_{s}}$$
	and 
	$$\dfrac{ng}{\mathrm{gcd}(\gamma_{\mathbf{x},I},ng)} = n_{1}gp_{1}^{j_{1}-m_{1}}\cdots p_{s}^{j_{s}-m_{s}},$$
	which implies
	$$\mathrm{rad}(\dfrac{ng}{\mathrm{gcd}(\gamma_{\mathbf{x},I},ng)}) = \mathrm{rad}(g)p_{1}^{z_{1}}\cdots p_{s}^{z_{s}},$$
	where $z_{k} = \mathrm{min}\{j_{k}-m_{k},1\}$ for any $1 \leq k \leq s$. Setting $\omega_{\mathfrak{m}}$ as above, then by Corollary \ref{coro 4}
	$$c_{ng}(\gamma_{\mathbf{x},I}) = \bigsqcup_{k=0}^{\omega_{\mathfrak{m}}-1}c_{ng/q^{\omega_{\mathfrak{m}}}}(\gamma_{\mathbf{x},I} q^{k})$$
	is the coarsest multiple equal-difference representation of $c_{ng}(\gamma_{\mathbf{x},I})$, where $c_{ng/q^{\omega_{\mathfrak{m}}}}(\gamma_{\mathbf{x},I} q^{k})$ denotes the $q^{\omega_{\mathfrak{m}}}$-cyclotomic coset modulo $ng$ containing $\gamma_{\mathbf{x},I} q^{k}$. And all these cosets $c_{ng/q^{\omega_{\mathfrak{m}}}}(\gamma_{\mathbf{x},I} q^{k})$ have the same common difference, given by $d_{\mathfrak{m}} = \frac{n\omega_{\mathfrak{m}}}{\tau_{s}(\mathfrak{m})}$.
	It then follows from Lemma \ref{lem 7} that $M_{h\gamma_{\mathbf{x},I}}(X)$ can be expressed as
	$$M_{h\gamma_{\mathbf{x},I}}(X) = \sum_{k=0}^{\omega_{\mathfrak{m}}}(-1)^{\omega_{\mathfrak{m}}-k}\sigma_{\omega_{\mathfrak{m}}-k}(\zeta_{d_{\mathfrak{m}}}^{h\gamma_{\mathbf{x},I}},\zeta_{d_{\mathfrak{m}}}^{h\gamma_{\mathbf{x},I} q},\cdots,\zeta_{d_{\mathfrak{m}}}^{h\gamma_{\mathbf{x},I} q^{\omega_{\mathfrak{m}}-1}})X^{\frac{\tau_{s}(\mathfrak{m}) k}{\omega_{\mathfrak{m}}}}.$$
\end{proof}

Now we are ready to state the main theorem in this section, which determines all $\lambda$-constacyclic codes of length $N$ over $\mathbb{F}_{q}$ by giving their generator polynomials.

\begin{theorem}
	The $\lambda$-constacyclic codes of length $N$ over $\mathbb{F}_{q}$ are
	$$\mathcal{C}_{\Omega} = (\prod_{\mathbf{x} \in [0,g_{1}^{y_{1}}-1]\times\cdots[0,g_{u}^{y_{u}}-1]}\prod_{I \in \mathcal{I}_{ng/n_{1}g}(\gamma_{\mathbf{x}})}M_{h\gamma_{\mathbf{x},I}}(X)^{\Omega_{\mathbf{x},I}}) \subseteq \mathcal{R}_{N,\lambda},$$
	where for any $\mathbf{x} \in [0,g_{1}^{y_{1}}-1]\times\cdots[0,g_{u}^{y_{u}}-1]$ and $I = (\mathfrak{m},\mathfrak{a},\mathfrak{t}) \in \mathcal{I}_{ng/n_{1}g}(\gamma_{\mathbf{x}})$ the exponent $\Omega_{\mathbf{x},I}$ varies the integer interval $[0,p^{v}]$.
\end{theorem}

Finally we conclude this section by illustrating the procedure given in Section \ref{sec 1} and \ref{sec 6} with a concrete example.

\begin{examples}
	\label{ex:91}
	Let
	\[
	q=5,\quad N=210,\quad \lambda=4\in\mathbb F_5^*.
	\]
	We determine all the $\lambda$-constacyclic codes of length $N$ over $\mathbb{F}_{q}$, along with their dimensions. Since
	\[
	\operatorname{ord}_{\mathbb F_5^*}(4)=2,
	\]
	we have
	\[
	N=5\cdot 42,\qquad
	n=42,\qquad
	g=\operatorname{ord}_{\mathbb F_5^*}(\lambda)=2.
	\]
	Moreover,
	\[
	\lambda'= \lambda^{g}=4^2=1.
	\]
	Hence
	\[
	X^{210}-4
	=
	\left(X^{42}-4\right)^5.
	\]
	Let $\zeta$ be a primitive $84$-th root of unity satisfying
	\[
	\zeta^{42}=4=-1\in\mathbb F_5.
	\]
	Then the defining set of $X^{42}-4$ is
	\[
	1+2\mathbb Z/84\mathbb Z,
	\]
	that is, the set of odd residue classes modulo $84$.
	
	Since $84=2^2\cdot21$, we first determine the $5$-cyclotomic cosets modulo the odd integer $21$ and then lift them from modulus $21$ to modulus $84$.
	
	It is trivial to check that the $5$-cyclotomic cosets modulo $21$ are
	\[
	\begin{aligned}
		c_{21}(0)&=\{0\},\\
		c_{21}(1)&=\{1,5,4,20,16,17\},\\
		c_{21}(2)&=\{2,10,8,19,11,13\},\\
		c_{21}(3)&=\{3,15,12,18,6,9\},\\
		c_{21}(7)&=\{7,14\}.
	\end{aligned}
	\]
	Thus we may take
	\[
	\mathcal{CR}_{21}=\{0,1,2,3,7\},
	\]
	with corresponding sizes
	\[
	|c_{21}(0)|=1,\quad
	|c_{21}(1)|=|c_{21}(2)|=|c_{21}(3)|=6,
	\quad
	|c_{21}(7)|=2.
	\]

	For $\gamma\in \mathcal{CR}_{21}$, let
	\[
	\tau=|c_{21}(\gamma)|.
	\]
	Since
	\[
	5^1\equiv1\pmod4,\quad
	5^6\equiv1\pmod4,\quad
	5^2\equiv1\pmod4,
	\]
	all five base cosets belong to the case
	\[
	5^\tau\equiv1\pmod4
	\]
	of Theorem~7.1. The corresponding values of
	\[
	v=v_2(5^\tau-1)
	\]
	are
	\[
	v=
	\begin{cases}
		2,&\gamma=0,\\
		3,&\gamma=1,2,3,7.
	\end{cases}
	\]
	
	For each $\gamma$, the cyclotomic $2$-adic integer is
	\[
	\phi(\gamma)=-\frac{\gamma}{21}\in\mathbb Z_2.
	\]
	Since we only need to lift from $21$ to $84=4\cdot21$, it is enough to determine
	the first two $2$-adic digits
	\[
	\phi(\gamma)_{\leq 1}
	=
	\varphi_0+2\varphi_1.
	\]
	They are determined by
	\[
	\gamma+21(\varphi_0+2\varphi_1)\equiv0\pmod4.
	\]
	As $21\equiv1\pmod4$, we obtain
	\[
	\begin{array}{c|c|c}
		\gamma & (\varphi_0,\varphi_1) &
		\phi(\gamma)_{\leq1}\\
		\hline
		0&(0,0)&0\\
		1&(1,1)&3\\
		2&(0,1)&2\\
		3&(1,0)&1\\
		7&(1,0)&1
	\end{array}
	\]
	
	By Theorem \ref{thm 5}, the principal sequence gives the coset
	\[
	c_{84}\bigl(\gamma+21\phi(\gamma)_{\leq1}\bigr),
	\]
	while the other lifts are obtained from
	\[
	c_{84}\bigl(\gamma+
	21(U(\gamma;m)+2^{m+1}t)\bigr),
	\]
	where
	\[
	U(\gamma;m)
	=
	\sum_{k=0}^{m-1}\varphi_k2^k
	+
	(1-\varphi_m)2^m.
	\]
	
	For $\gamma=0$, we have
	\[
	\phi(0)_{\leq1}=0,
	\quad
	U(0;0)=1,
	\quad
	U(0;1)=2.
	\]
	Hence the degree-$2$ lifts are
	\[
	c_{84}(21),\quad
	c_{84}(63).
	\]

	For $\gamma=1$, we have
	\[
	\phi(1)_{\leq1}=3,
	\quad
	U(1;0)=0,
	\qquad
	U(1;1)=1.
	\]
	Thus the four lifts are 
	\[
	c_{84}(1),\quad c_{84}(43).
	\]
	
	For $\gamma=2$, we have
	\[
	\phi(2)_{\leq1}=2,
	\quad
	U(2;0)=1,
	\qquad
	U(2;1)=0.
	\]
	The lifts are 
	\[
	c_{84}(11),\quad c_{84}(13).
	\]
	
	For $\gamma=3$, we have
	\[
	\phi(3)_{\leq1}=1,
	\quad
	U(3;0)=0,
	\qquad
	U(3;1)=3.
	\]
	The lifts are 
	\[
    c_{84}(3),\quad c_{84}(9).
	\]
	
	Finally, for $\gamma=7$, we have
	\[
	\phi(7)_{\leq1}=1,
	\quad
	U(7;0)=0,
	\qquad
	U(7;1)=3.
	\]
	The lifts are 
	\[
	c_{84}(7),\quad c_{84}(49).
	\]
	
	Consequently, the complete set of $5$-cyclotomic cosets modulo $84$
	contained in the defining set $1+2\mathbb Z/84\mathbb Z$ is
	\[
		\begin{aligned}
			c_{84}(21)&=\{21\},\\
			c_{84}(63)&=\{63\},\\
			c_{84}(7)&=\{7,35\},\\
			c_{84}(49)&=\{49,77\},\\
			c_{84}(1)&=\{1,5,25,41,37,17\},\\
			c_{84}(3)&=\{3,15,75,39,27,51\},\\
			c_{84}(9)&=\{9,45,57,33,81,69\},\\
			c_{84}(11)&=\{11,55,23,31,71,19\},\\
			c_{84}(13)&=\{13,65,73,29,61,53\},\\
			c_{84}(43)&=\{43,47,67,83,79,59\}.
	\end{aligned}
	\]
	Thus a complete set of representatives is
	\[
	\mathcal{R}=\{21,63,7,49,1,3,9,11,13,43\}.
	\]
	Their sizes are
	\[
	1,1,2,2,6,6,6,6,6,6,
	\]
	respectively, and hence
	\[
	1+1+2+2+6+6+6+6+6+6=42.
	\]
	This agrees with the fact that $X^{42}-4$ has degree $42$.
	
	For $\gamma\in \mathcal{R}$, define
	\[
	M_\gamma(X)
	=
	\prod_{i=0}^{|c_{84}(\gamma)|-1}
	(X-\zeta^{5^i\gamma}).
	\]
	The two linear factors are
	\[
	M_{21}(X)=X-2,
	\qquad
	M_{63}(X)=X-3,
	\]
	and the two quadratic factors are
	\[
	M_7(X)=X^2+2X+4,
	\qquad
	M_{49}(X)=X^2-2X+4.
	\]
	The factors of degree $6$ are
	\[
	\begin{aligned}
		M_1(X)
		&=X^6-2X^4+X^3+2X^2-1,\\
		M_3(X)
		&=X^6-2X^5-X^4+2X^3+X^2-2X-1,\\
		M_9(X)
		&=X^6+2X^5-X^4-2X^3+X^2+2X-1,\\
		M_{11}(X)
		&=X^6+2X^5+2X^4+2X^3-2X^2+2X-1,\\
		M_{13}(X)
		&=X^6-2X^5+2X^4-2X^3-2X^2-2X-1,\\
		M_{43}(X)
		&=X^6-2X^4-X^3+2X^2-1.
	\end{aligned}
	\]
	Therefore
	\[
	\begin{aligned}
		X^{42}-4
		={}&(X-2)(X-3)
		(X^2+2X+4)(X^2-2X+4)\\
		&\cdot M_1(X)M_3(X)M_9(X)
		M_{11}(X)M_{13}(X)M_{43}(X).
	\end{aligned}
	\]
	
	Since
	\[
	X^{210}-4=(X^{42}-4)^5,
	\]
	every monic divisor of $X^{210}-4$ is of the form
	\[
	g(X)=
	\prod_{\gamma\in R}M_\gamma(X)^{e_\gamma},
	\qquad
	0\leq e_\gamma\leq5.
	\]
	Hence every ideal
	\[
	\mathcal{C}=(g(X))
	\subseteq
	\mathbb F_5[X]/(X^{210}-4)
	\]
	is a $4$-constacyclic code of length $210$, and
	\[
	\dim_{\mathbb F_5}\mathcal{C}
	=
	210-\deg g(X)
	=
	210-\sum_{\gamma\in R}
	e_\gamma |c_{84}(\gamma)|.
	\]
	In particular, there are
	\[
	6^{10}=60466176
	\]
	distinct $4$-constacyclic codes of length $210$ over $\mathbb F_5$.
\end{examples}

\appendix
\renewcommand{\thesection}{Appendix~\Alph{section}}

\section{The multiple equal-difference representations of cyclotomic cosets}
In this appendix we recall the definition and some results of the multiple equal-difference representations of cyclotomic cosets. For a complete treatment the readers are referred to \cite{Zhu_2}. Since in this appendix we will deal with $q^{t}$-cyclotomic coset for some positive integer $t$, we write $c_{n/q^{t}}(\gamma)$ for the $q^{t}$-cyclotomic coset modulo $n$ containing $\gamma$.

\begin{definition}
	Let $E$ be a subset of $\mathbb{Z}/n\mathbb{Z}$ with order $| E | =\tau_{E}$. We say that $E$ is of equal difference if $\tau_{E} \mid n$ and $E$ can be presented as 
	$$E = \{\gamma,\gamma+\frac{n}{\tau_{E}},\cdots,\gamma+(\tau_{E}-1)\frac{n}{\tau_{E}}\}.$$
	The integer $\frac{n}{\tau_{E}}$ is called the common difference of $E$. In particular, an equal-difference $q$-cyclotomic coset modulo $n$ is a $q$-cyclotomic coset modulo $n$ that is of equal difference.
\end{definition}

For any $\gamma \in \mathbb{Z}/n\mathbb{Z}$ we set $n_{\gamma} = \frac{n}{\mathrm{gcd}(\gamma,n)}$. It is trivial to see that for any $\gamma$ and $\gamma^{\prime}$ lying in the same $q$-cyclotomic coset modulo $n$ it holds that $n_{\gamma} = n_{\gamma^{\prime}}$. The criterion for a $q$-cyclotomic coset modulo $n$ being of equal difference is given below.

\begin{proposition}
	The cyclotomic coset $c_{n/q}(\gamma)$ is of equal difference if and only if the following two conditions are satisfied:
	\begin{description}
		\item[(\romannumeral1)] $\mathrm{rad}(n_{\gamma}) \mid q-1$;
		\item[(\romannumeral2)] $q \equiv 1 \pmod{4}$ if $8 \mid n_{\gamma}$.
	\end{description}
\end{proposition}

In general, a cyclotomic coset $c_{n/q}(\gamma)$ is not necessarily of equal difference, however it turns out that $c_{n/q}(\gamma)$ can always be expressed as a disjoint union of equal-difference subsets, which is called an equal-difference decomposition of $c_{n/q}(\gamma)$. If, moreover, the decomposition is in the form
$$c_{n/q}(\gamma) = \bigsqcup_{j=0}^{\mathrm{gcd}(t,\tau)-1}c_{n/q^{t}}(\gamma q^{j})$$
for some $t \in \mathbb{N}^{+}$, then it is called a multiple equal-difference representation of $c_{n/q}(\gamma)$.

Let $c_{n/q}(\gamma)$ be a $q$-cyclotomic coset modulo $n$, with $n_{\gamma} = \frac{n}{\mathrm{gcd}(\gamma,n)}$. Define a positive integer $\omega_{\gamma}$ by
\begin{equation*}
	\omega_{\gamma} = \left\{
	\begin{array}{lcl}
		2\mathrm{ord}_{\mathrm{rad}(n_{\gamma})}(q), \quad \mathrm{if} \ q^{\mathrm{ord}_{\mathrm{rad}(n_{\gamma})}(q)} \equiv 3 \pmod{4} \ \mathrm{and} \ 8 \mid n_{\gamma};\\
		\mathrm{ord}_{\mathrm{rad}(n_{\gamma})}(q), \quad \mathrm{otherwise}.
	\end{array} \right.
\end{equation*}
The following theorem determines all the multiple equal-difference representations of $c_{n/q}(\gamma)$.

\begin{theorem}
	All the distinct multiple equal-difference representations of $c_{n/q}(\gamma)$ are
	$$c_{n/q}(\gamma) = \bigsqcup_{j=0}^{t-1}c_{n/q^{t}}(\gamma q^{j}),$$
	where $t$ ranges over all positive divisors of $\tau = |c_{n/q}(\gamma)|$ that are divided by $\omega_{\gamma}$.
\end{theorem}

\begin{corollary}\label{coro 4}
	The coarsest multiple equal-difference representation of $c_{n/q}(\gamma)$ is
	$$c_{n/q}(\gamma) = \bigsqcup_{j=0}^{\omega_{\gamma}-1}c_{n/q^{\omega_{\gamma}}}(\gamma q^{j}),$$
	while the finest multiple equal-difference representation of $c_{n/q}(\gamma)$ is the trivial decomposition
	$$c_{n/q}(\gamma) = \bigsqcup_{j=0}^{\tau-1}\{\gamma q^{j}\}.$$
\end{corollary}

The multiple equal-difference representations of $q$-cyclotomic cosets modulo $n$ are particularly of interest to us as they encodes the irreducible factorizations of $X^{n}-1$ into binomials over finite extensions of $\mathbb{F}_{q}$. The explicit statements are given below.

\begin{theorem}
	\begin{description}
		\item[(1)] A $q$-cyclotomic coset $c_{n/q}(\gamma)$ modulo $n$ is of equal-difference if and only if the induced polynomial
		$$M_{\gamma}(X) = (X-\zeta_{n}^{\gamma})(X-\zeta_{n}^{\gamma q})\cdots(X-\zeta_{n}^{\gamma q^{\tau-1}})$$
		is a binomial.
		\item[(2)] For a positive integer $t$, the following statements are equivalent:
		\begin{description}
			\item[(\romannumeral1)] $\omega_{\gamma} \mid t$;
			\item[(\romannumeral2)] $c_{n/q}(\gamma) = \bigsqcup\limits_{j=0}^{\mathrm{gcd}(t,\tau)-1}c_{n/q^{t}}(\gamma q^{j})$ is a multiple equal-difference representation of $c_{n/q}(\gamma)$;
			\item[(\romannumeral3)] $M_{\gamma}(X)$ is factorized into irreducible binomials over $\mathbb{F}_{q^{t}}$.
		\end{description}
		\item[(3)] For a positive integer $t$, the following statements are equivalent:
		\begin{description}
			\item[(\romannumeral1)] $\omega_{1} \mid t$;
			\item[(\romannumeral2)] $c_{n/q}(\gamma) = \bigsqcup\limits_{j=0}^{\mathrm{gcd}(t,\tau)-1}c_{n/q^{t}}(\gamma q^{j})$ is a multiple equal-difference representation of $c_{n/q}(\gamma)$ for every $q$-cyclotomic coset $c_{n/q}(\gamma)$ modulo $n$;
			\item[(\romannumeral3)] $X^{n}-1$ is factorized into irreducible binomials over $\mathbb{F}_{q^{t}}$.
		\end{description}
	\end{description}
\end{theorem}

\begin{corollary}
	The field $\mathbb{F}_{q^{\omega_{1}}}$ is the smallest extension of $\mathbb{F}_{q}$ over which $X^{n}-1$ is factorized into irreducible binomials.
\end{corollary}
	
Let $c_{n/q}(\gamma)$ be a $q$-cyclotomic coset modulo $n$ with size $\tau$, and let $n_{\gamma}$ and $\omega_{\gamma}$ be defined as above. By Corollary \ref{coro 4}
$$c_{n/q}(\gamma) = \bigsqcup_{j=0}^{\omega_{\gamma}-1}c_{n/q^{\omega_{\gamma}}}(\gamma q^{j})$$
is the coarsest multiple equal-difference representation of $c_{n/q}(\gamma)$. For each $0 \leq j \leq \omega_{\gamma}-1$ the cosets $c_{n/q^{\omega_{\gamma}}}(\gamma q^{j})$ has size $\frac{\tau}{\omega_{\gamma}}$, thus it is an equal-difference set with common difference $d = \frac{n\omega_{\gamma}}{\tau}$. It follows that the induced irreducible polynomial $\widetilde{M}_{\gamma q^{j}}(X)$ over $\mathbb{F}_{q^{\omega_{\gamma}}}$ is 
$$\widetilde{M}_{\gamma q^{j}}(X) = (X-\zeta_{n}^{\gamma q^{j}})(X-\zeta_{n}^{\gamma q^{j}+d})\cdots (X-\zeta_{n}^{\gamma q^{j}+(\frac{n}{d}-1)d}) = X^{\frac{n}{d}} - \zeta_{d}^{\gamma q^{j}}.$$
Hence we have 
$$M_{\gamma}(X) = \prod_{j=0}^{\omega_{\gamma}-1}\widetilde{M}_{\gamma q^{j}}(X) = \sum_{j=0}^{\omega_{\gamma}}(-1)^{\omega_{\gamma}-j}\sigma_{\omega_{\gamma}-j}(\zeta_{d}^{\gamma},\zeta_{d}^{\gamma q},\cdots,\zeta_{d}^{\gamma q^{\omega_{\gamma}-1}})X^{\frac{\tau j}{\omega_{\gamma}}},$$
where $\sigma_{\omega_{\gamma}-j}$ denotes the $(\omega_{\gamma}-1)$-th elementary symmetric polynomial. This proves the lemma below.

\begin{lemma}\label{lem 7}
	The irreducible polynomial $M_{\gamma}(X)$ induced by $c_{n/q}(\gamma)$ can be expressed as
	$$M_{\gamma}(X) = \sum_{j=0}^{\omega_{\gamma}}(-1)^{\omega_{\gamma}-j}\sigma_{\omega_{\gamma}-j}(\zeta_{d}^{\gamma},\zeta_{d}^{\gamma q},\cdots,\zeta_{d}^{\gamma q^{\omega_{\gamma}-1}})X^{\frac{\tau j}{\omega_{\gamma}}},$$
	where $\sigma_{\omega_{\gamma}-j}(\zeta_{d}^{\gamma},\zeta_{d}^{\gamma q},\cdots,\zeta_{d}^{\gamma q^{\omega_{\gamma}-1}})$ is the evaluation of the $(\omega_{\gamma}-1)$-th elementary symmetric polynomial at the point $(\zeta_{d}^{\gamma},\zeta_{d}^{\gamma q},\cdots,\zeta_{d}^{\gamma q^{\omega_{\gamma}-1}})$.
\end{lemma}
	
\section*{Acknowledgment}
The first author was supported by Basic Research Program Young Scientists Guidance Project of Guizhou Province (QN[2025]186). 
The authors are grateful to the editor and the anonymous referees for taking time to read and 
comment on the article very carefully and insightfully. Their nice comments and suggestions helped 
us to improve the article very much.

\section*{Data availability}
Data sharing not applicable to this article as no datasets were generated or analysed during the current study.

\section*{Declaration of competing interest}
The authors declare that we have no known competing financial interests or personal relationships 
that could be perceived to influence the work reported in this paper.

\end{document}